\documentclass[a4paper,11pt]{article}

\usepackage{colortbl}
\usepackage{multirow}%
\usepackage{fullpage}
\usepackage{amsmath}%
\usepackage{amsthm}%
\usepackage{amsfonts}%
\usepackage{amssymb}%
\usepackage{graphicx}
\usepackage{enumerate}
\usepackage[numbers]{natbib}
\usepackage[colorlinks=true,linkcolor=black,citecolor=blue]{hyperref}
\usepackage[titletoc,title]{appendix}
\usepackage{subfig}
\usepackage{changes}

\newtheorem{theorem}{Theorem}

\newtheorem{proposition}[theorem]{Proposition}

\newtheorem{example}{Example}
\newtheorem{app_theorem}{Theorem}[section]
\newtheorem{app_lemma}[app_theorem]{Lemma}

\newcommand{\E}{\mathbb{E}}

\newcommand{\diff}{\mathrm{d}}
\newcommand{\fthetahat}{f_{\hat{\theta}_n}}
\newcommand{\alphahat}{\hat{\alpha}_n}

\DeclareMathOperator{\var}{var}

\DeclareMathOperator{\argmax}{argmax}
\DeclareMathOperator{\np}{(np)}
\DeclareMathOperator{\p}{(p)}

\newcounter{ass:H}
\newcounter{ass:A}
\newcounter{ass:B}

\usepackage{setspace}

\begin{document}

{
  \title{Parametric versus nonparametric: the fitness coefficient}
\vspace{2cm}
  \author{Gildas Mazo\thanks{The research by G. Mazo was partially funded by a ``Projet de Recherche" of the ``Fonds de la Recherche Scientifique --- FNRS" (Belgium).} \\
    \small MaIAGE, INRA, Universit\'e Paris-Saclay\\
    Fran\c{c}ois Portier\thanks{
    The author gratefully acknowledges support from the Fonds de la Recherche Scientifique (FNRS) A4/5 FC 2779/2014-2017 No.\ 22342320, from the contract ``Projet d'Act\-ions de Re\-cher\-che Concert\'ees'' No.\ 12/17-045 of the ``Communaut\'e fran\c{c}aise de Belgique''.}\\
    \small LTCI, T\'el\'ecom ParisTech, Universit\'e Paris-Saclay}
  \maketitle
}





\begin{abstract}
The fitness coefficient, introduced in this paper, results from a competition between  parametric and nonparametric density estimators within the likelihood of the data. As illustrated on several real datasets, the fitness coefficient generally agrees with p-values but  is easier to compute and interpret. Namely, the fitness coefficient can be interpreted as the proportion of data coming from the parametric model. Moreover, the fitness coefficient can be used to build a semiparamteric compromise which improves inference over the parametric and nonparametric approaches.
From a theoretical perspective, the fitness coefficient is shown to converge in probability to one if the model is true and to zero if the model is false. From a practical perspective, the utility of the fitness coefficient is illustrated on real and simulated datasets.

\end{abstract}





\section{Introduction}
A challenge of data analysis is to assess the quality of a model. The traditional approach relies on goodness-of-fit tests where, loosely speaking, the ability of a model to fit the data is measured through distances between the observed values and the values expected under the model.
Examples include the classical Pearson's chi-squared test~\cite{akritas:1988}, the Kolmogorov or Cramer-von-Mises goodness-of-fit  tests~\cite{anderson+h+t:1994,dagostinoStephensBook,durbin:1973}, or
likelihood-ratio based statistics~\cite{chen+vk:2009,claeskens+h:2004,vuong:1989} (see~\cite{chen+vk:2009} for an empirical likelihood approach).

In this context, p-values have emerged as natural instruments to measure the amount of evidence in favor of the model. However, the use of p-values is subjected to several difficulties:
(i) their calculation might require computationally intensive strategies as the bootstrap~\cite{MR3647588,MR2968394,stute1993bootstrap}; (ii)~interpretation is notoriously  difficult~\citep{murtaugh:2014,schervish:1996} as emphasized again in a recent ASA statement~\citep{wasserstein:2016}
; (iii) whenever some evidence has been found against the model, no information is delivered to improve inference.

In this paper, we introduce \emph{the fitness coefficient}, a new criterion for simultaneously measuring the amount of evidence of a model and improving inference. Our goal is to provide an alternative approach to the use of p-values in goodness-of-fit testing that is no longer sensitive to the difficulties (i)-(iii). 

Let $X_1,\dots,X_n$ be independent $d$-variate observations with common density $f_0$. Let $\mathcal{P} = \{f_\theta\,:\,\theta\in \Theta\}$ be a family of probability density functions representing the model.  Given the maximum likelihood estimator $\fthetahat$ (based on the model), and the standard kernel density estimator $\hat f_{n}$ (free from the model) 
with kernel $K: \mathbb R^d \to \mathbb R_+$ and bandwidth $h_n>0$, define the fitness coefficient $\alphahat$ as 
\begin{align}\label{eq:our method}
  \alphahat \in \underset{\alpha \in [0,1]}{\argmax}\ 
  \sum_{i=1}^n \log\left( \alpha \fthetahat(X_i) + (1-\alpha)\hat f_{i,n}^{\text{LR}}   \right),
\end{align}
where $\hat f_{i,n}^{\text{LR}} $ is called the \textit{leave-and-repair} (LR) kernel estimate of $f_0(X_i)$ and is given by
\begin{align}\label{def:loo}
\hat f_{i,n}^{\text{LR}} = \left( \frac{1}{(n-1)h_n^d} \sum_{j\neq i}   K\left( \frac {X_i-X_j}{h_n} \right) \right) +  \Delta_n q(X_i) , 
\end{align}
with $\Delta_n\geq 0$ and $q: \mathbb R^d \to \mathbb R_+$. The LR estimate is a modification of the well-known \textit{leave-one-out} (LOO) estimate usually employed in cross-validation procedures \cite{hall:1987} and semiparametric estimation~\cite{delecroix+h+p:2006}. 

The fitness coefficient has the following advantages.

\paragraph{(i) The fitness coefficient is easy to compute.} It is the minimizer of a simple one dimensional concave function.

\paragraph{(ii) The fitness coefficient is a measure of model quality.}

As seen in~(\ref{eq:our method}), the fitness coefficient $\hat \alpha_n$ follows from a competition between the parametric and the nonparametric approach so as to maximize the likelihood of the observations. Hence, whenever the model is \textit{sufficiently} true, we expect a value of $\alphahat$ \textit{relatively} close to one. This is because the parametric estimator is likely to be more accurate than the nonparametric one. On the opposite, whenever the model is wrong, we expect a value of $\alphahat$ close to zero. 
Because $\alphahat \fthetahat + (1-\alphahat) \hat{f}_n$ is a mixture distribution between the parametric and the nonparametric estimates, the fitness coefficient is interpreted as the proportion of data distributed under the model. For instance, if one draws a bootstrap sample from the combination $\alphahat \fthetahat + (1-\alphahat) \hat{f}_n$ then the fitness coefficient $\alphahat$ is the proportion of data drawn from the fitted model $\fthetahat$. Therefore, the less the value of $\alphahat$, the less the bootstrap sample shall be ``contaminated'' by the nonparametric part of the combination. 

\begin{figure}
 \centering
 \includegraphics[width=0.5\textwidth]
 {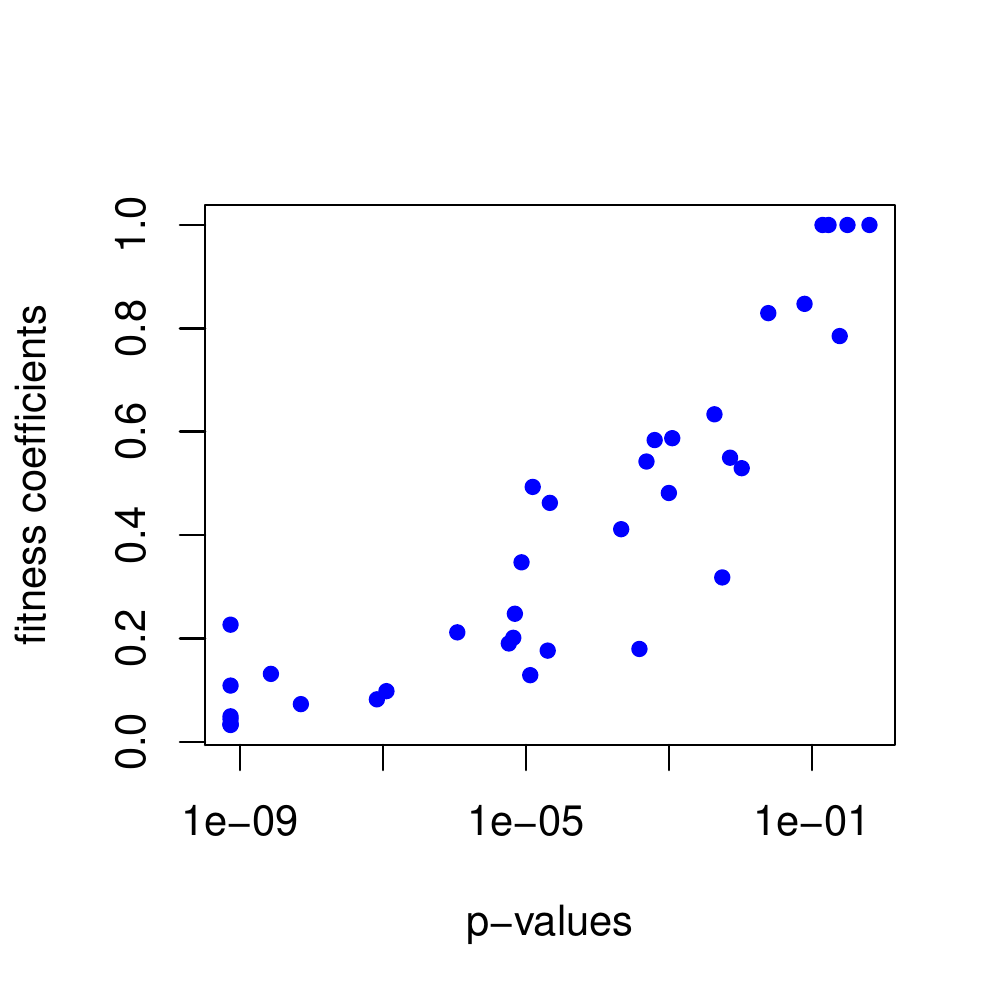}
 \caption{Estimated p-values (on a logarithmic scale) of the Cramer-Von-Mises goodness-of-fit normality test against the values of the fitness coefficient for the CAC40 data.}
 \label{fig:cac40}
\end{figure}

To show the capability of the fitness coefficient, we compare it to p-values on a real data example. Consider the problem of testing whether 
a given sample comes from a normal distribution. Specifically, we have
38 samples each consisting of $n=409$ financial returns of a company from the French stock market CAC40 and we wish to measure the quality of the normal model for each of these samples. 
On the one hand a goodness-of-fit test based on the Cramer-Von-Mises statistic~\cite{dagostinoStephensBook}
is carried out. On the other hand the fitness coefficient defined by~(\ref{eq:our method}) is computed with a Gaussian kernel $K$,  $\Delta_n=1/n$, $q(x)=t_\nu(x/100)$ where $t_\nu$ is the density of a Student-t distribution with $\nu=3$ degrees of freedom, and $h_n$ given by~\cite{silvermanbook1998} (p. 48).
In Figure~\ref{fig:cac40}, we plotted the values of the fitness coefficient against the p-values on a logarithmic scale.  We can see a clear positive dependence relationship: large values for the fitness coefficient correspond to large p-values. 
This suggests that, if one had used p-values to assess the fitness of the normal model, he or she could have done so with the fitness coefficient.

The quality criterion induced by the fitness coefficient is different than that of information criteria~\cite{burnham+a:2003,claeskens+h:2008} such as the Akaike information criterion~\cite{akaike1974new} or the Bayesian information criterion ~\cite{schwarz:1978} which focus on the \emph{relative} performances \emph{between} models. Note that convex parametric combinations recently have been proposed in the Bayesian literature~\cite{kamary2014testing} to assess the fitness of a certain parametric model against another.

\paragraph{(iii) The fitness coefficient is useful to get robust semiparametric estimators.}
The fitness coefficient offers a natural semiparametric alternative $\alphahat \fthetahat + (1-\alphahat) \hat{f}_n$ for estimating the probability density function $f_0$ of the observations. The idea of forming such a convex combination to get an estimator robust to misspecification while retaining a performance comparable to parametric estimators when the true density is close to the model was originally developed by Olkin and Spiegelman~\cite{olkinSpiegelman1987}. Their method, referred to as the \emph{OS method}, consists of computing
\begin{align}\label{eq:OSmethod}
  \hat{\alpha}_n^{\text{OS}} \in \underset{\alpha \in [0,1]}{\argmax}\ 
  \sum_{i=1}^n \log\left( \alpha \fthetahat(X_i) + (1-\alpha) \hat{f}_n(X_i) \right),
\end{align}
where $\hat{f}_n$ is the standard kernel density estimator. The OS method and the LR method given in~(\ref{eq:our method}) differ in the choice of the nonparametric estimator in the combination. The OS method was noticed to be sensitive to the choice of the bandwidth~\cite{faraway1990,rahmanEtAl-aNote-1997-brazilian}.

Rather than considering the likelihood of the observations, some authors \cite{leeSoleymani-mixingEstimators-2015-JASA,rahmanEtAl-aNote-1997-brazilian,leeSoleymani-localSemiParametric-2014-JSPI} investigate strategies based on the mean squared error between the combination $\alpha \fthetahat + (1-\alpha) \hat{f}_n$ and the true density $f_0$, but then the solution depends on the unknown distribution and hence heavy bootstrap methods need to be employed.

To improve inference, there exist also other approaches than that of forming a convex combination between the parametric and nonparametric estimators. The locally parametric nonparametric estimation is developed for instance in \cite{hjort+j:1996,hjort+m+vk:2017hybrid,talamakrouni+g+vk:2017}, but is less appealing from the point of view of model quality assessment because they do not provide any ``fitness coefficient''.

\bigskip

\paragraph{Main contributions.}

By introducing the fitness coefficient, we provide a new measure for assessing the quality of a model and an alternative to the OS method to get robust semiparametric estimators. Under mild conditions, the fitness coefficient $\alphahat$ is shown to converge in probability to one if the model is true and zero otherwise, a property called \emph{consistency}. Even if the fitness coefficient is maximizing some objective function (over $\alpha\in [0,1]$), classical results from M-estimation theory does not apply because, when the model is true, the limiting objective function is independent from $\alpha$. The proposed approach follows from a fine comparison between the rates of convergence of $\fthetahat$ and $ \hat{f}_n$. We moreover provide examples of densities $f_0$ that satisfy our set of assumptions. Using some real data as well as extensive simulations, we observed that the LR approach is more stable than the OS approach and leads to more accurate inference. This is in agreement with our theoretical analysis, which cannot include the OS method as an example.

\paragraph{Outline.}
 In Section~\ref{sec:def+framework}, we introduce some quantities of interest and motivate the use of the LR estimator $\hat f_{i,n}^{\text{LR}}$ to compute the fitness coefficient $\alphahat$. The consistency of the fitness coefficient is stated in Section~\ref{sec:main} where some examples are given. In Section~\ref{sec:numerical_results}, numerical experiments are designed to measure the robustness of the fitness coefficient and the performance of the corresponding density estimators.
All the proofs are postponed to the Appendix.

\section{The leave-and-repair estimator}\label{sec:def+framework}

The aim of this section is to define and motivate the use of the leave-and-repair (LR) estimator $\hat f_n^{\text{LR}}$~(\ref{def:loo}) introduced in the definition of the fitness coefficient $\hat \alpha_n$~(\ref{eq:our method}). 
Compared with the OS method given in (\ref{eq:OSmethod}), the use of the LR estimator might seem irrelevant at first view but in fact plays an important role to ensure the good behavior of the fitness coefficient. 

The kernel density estimator of $f_0$ at $x\in \mathbb R^d$ is given by
\begin{align*}
\hat f_{n}(x) =\frac{1}{nh_n^d} \sum_{i=1}^n   K\left( \frac {x-X_i}{h_n} \right).
\end{align*}
For any $h>0$, define the function $f_h$ as the convolution product between $K_h(\cdot) = K(\cdot /h) /h^d$ and $f_0$, that is, $f_{h} (x)  = (K_h\star f_0) (x)$, $x\in \mathbb R^d$. Note that $f_{h_n} (x)$ is the expected value of $\hat f_n (x)$. In other words, $f_{h_n} (x) = \mathbb E [ \hat f_n (x)] $. But since $\mathbb E [ \hat f_n(X_i) \vert X_i ] = f_{h_n}(X_i) + K(0)/(nh_n^d)$, we see that $\hat f_n (X_i)$ has a positive $h_n$-dependent bias when estimating $f_{h_n}(X_i)$, conditionally on $X_i$. When studying the estimator decomposition, this bias term spreads to the diagonal terms of some $U$-statistics and gives rise, in the end, to some non-negligible terms. This phenomenon is common in semiparametric statisitcs, and has been noticed for instance in Remark 4 in \cite{portier+s:2015}.

To overcome the undesirable effects caused by this bias term, the \textit{leave-one-out} (LOO) estimator of $f_{h_n}(X_i)$, given by
\begin{align*}
\hat{f}_{n,i}^{\text{LOO}} = \frac{1}{(n-1)h_n^d} \sum_{j\neq i}   K\left( \frac {X_i-X_j}{h_n} \right),
\end{align*}
has been successfully used in several cross-validation procedures aiming at selecting the bandwidth, either based on the likelihood~\cite{habbema+h+v:1974,marron:1985,hall:1987} or on the mean squared error \cite{rudemo:1982,stone:1984} (see \cite{marron:1987} for a comparison). Since then, LOO estimators have been frequently used in semiparametric studies \cite{delecroix+h+p:2006}. 

The LR estimator proposed in this paper is inspired, but different, from the LOO estimator. In view of~(\ref{def:loo}), the LR estimator satisfies 
\begin{align*}
\hat f_{n,i}^{\text{LR}} = \hat{f}_{n,i}^{\text{LOO}} +  \Delta_n q(X_i).
\end{align*}
If $\Delta_n =0$ the LR estimator is equal to the LOO estimator. If $q= K(0)$ and $\Delta_n  = 1/((n-1)h_n^d)$ the LR estimator is equal to $(n/(n-1))\hat f_n(X_i)$. In general, one can think of $\Delta_n$ as of order $1/n$ and $q$ as a density, yielding that $q(X_1)>0$ has probability $1$.

The heuristic for using the LR estimator $\hat f_{n,i}^{\text{LR}}$ instead of the LOO estimator $\hat f_{n,i}^{\text{LOO}}$ is as follows. It is well-known that the Kullback-Leibler divergence of kernel density estimates depends crucially on the tails of the true distribution $f_0$~\citep{hall:1987,schuster+g:1981}. As shown in~\cite{hall:1987}, if the tail is too heavy and the kernel $K(x)$ vanishes too quickly as $x$ becomes large then the Kullback-Leibler divergence associated to the kernel density estimate goes to minus infinity. This is because some of the $\hat f_{n,i} ^{\text{LOO}}$, $i=1,\ldots,n$, might have very small values (possibly zero), leading to very large values (possibly infinite) for some of the $\log(\hat f_{n,i} ^{\text{LOO}})$. These values are involved in the computation of the Kullback-Liebler divergence and play an important role in our proofs when dealing with the likelihood of the nonparametric estimate. We built the LR estimator to overcome this issue by simply adding $\Delta_n q(X_i)$ to the LOO estimator. We coined the term \textit{leave-and-repair} because the term $\Delta_n q(X_i)$ \emph{repairs} the LOO estimator. Since $\hat f _{n,i} ^{\text{LR}} \ge \Delta_n q(X_i)$, the LR estimator is not subjected to the difficulties of the LOO estimator. By adding the term $\Delta_n q(X_i)$ in~(\ref{def:loo}), however, a biais is introduced: now one has $E [ \hat f _{n,i} ^{\text{LR}} \vert X_i ] - f_{h_n}(X_i) = \Delta_n q(X_i)$. Thus, there is a biais-variance tradeoff controlled by the sequence $\Delta_n$ that must go to zero slowly enough to keep $\hat f _{n,i} ^{\text{LR}}$ away from zero but also fast enough to keep the biais as small as possible. The right compromise is given in the conditions in Theorem~\ref{main:result} (for instance $\Delta_n=1/n$ is one possibility).

Concerning the parametric estimator $\fthetahat$, we follow \cite{olkinSpiegelman1987} by considering the \textit{maximum likelihood estimator}. Let 
$\mathcal P = \{ f_\theta\ : \ \theta \in \Theta \}$ be the parametric model where $\Theta \in \mathbb R^p$ is such that for each $\theta\in \Theta$, $f_\theta : \mathbb R^d\to \mathbb R^+  $ is a measurable function satisfying $\int f_\theta(x) \,\diff x =1$. The maximum likelihood estimator of $f_0$ based on $\mathcal P$ and $X_1,\ldots, X_n$ is $f_{\hat \theta_n}$ where $\hat \theta_n$ (when it exists; this is further assumed) is defined as
\begin{align*}
\hat \theta_n \in \underset{\theta\in \Theta}{\argmax} \, \sum_{ i=1 }^n \log( f_\theta(X_i) ) .
\end{align*}
The good behaviour of the maximum likelihood estimator is subjected to the assumption that $ f_0\in \mathcal P$, that is, there exists $\theta_0\in \Theta $ such that $f_0 = f_{\theta_0}$.

To conclude the section, we consider existence and uniqueness of the fitness coefficient $\hat \alpha_n$. The existence follows from the use of the LR estimator $\hat f_{n,i}^{\text{LR}}$. Uniqueness of $\hat \alpha_n$ is obtained under the mild requirement that the parametric and nonparametric estimators are distinguishable on the observed data.

\begin{proposition}\label{prop:existence}
Suppose that $q(X_1)>0$ a.s. and $f_{\hat \theta_n}(X_i)\neq \hat f _{n,i} ^{\text{LR}}$ for at least one $i\in\{1,\ldots,n\}$. Then the fitness coefficient exists and is unique.
\end{proposition}

The proof is given in Appendix \ref{app:proof_prop}.

\section{Consistency of the fitness coefficient}\label{sec:main}

Recall that consistency of the fitness coefficient $\hat \alpha_n$ means $\alphahat \to 1$ in probability if $f_0 \in \mathcal P$ and $\alphahat \to 0$ if $f_0 \notin \mathcal P$, where $f_0$ is the true underlying density and 
$\mathcal{P}$ is the parametric model.
Section~\ref{sec:subsec:main-result} and Section~\ref{sec:subsec:examples-H3} contain the main consistency theorem and some examples satisfying our set of assumptions, respectively.

\subsection{Assumptions and main result}\label{sec:subsec:main-result}
Let $\|\cdot\|_2$ be the Euclidean norm and for any set $S\subset \mathbb R^d$ and any function $f:S \to \mathbb R$, define the sup-norm as $\|f\|_S = \sup_{x\in S}|f(x)|$. Denote by $\lambda$ the Lebesgue measure on $\mathbb{R}^d$. Introduce the density level sets
\begin{align*}
S_t = \{x\in \mathbb R^d \,:\,  f_0(x) >t \},\qquad t\geq 0.
\end{align*}
We shall assume the following.

\begin{enumerate}[(H1)]
\item\label{ass:base_density} The density $f_0$ is bounded and continuous on $\mathbb R^d$ and the gradient $\nabla f_0$ of $f_0$ is bounded on $\mathbb R^d$, and satisfies, for every $x\in \mathbb R^d$ and $u\in [-1,1]^d$,
\begin{align*}
|f_0(x+u) - f_0(x) - u^T \nabla f_0(x)| \leq \| u\|_2^2 g(x),
\end{align*}
where $g$ is positive, bounded, integrable and $\int  {g(x)^2 }/{f_0(x)} \,\diff x < \infty$.
\setcounter{ass:H}{\value{enumi}}
\end{enumerate}

\begin{enumerate}[(H1)]
\setcounter{enumi}{\value{ass:H}}
\item \label{ass:kernel}
The kernel function $K:\mathbb R^d \rightarrow \mathbb R^+$ integrates to $1$ and takes one of the two following forms,
 \begin{align*}
  (a) \quad K(x) \propto  K^{(0)} (\|x\|_2), \qquad \text{or} \qquad (b) \quad  K(x) \propto \prod_{k=1}^d  K^{(0)}(|x_k|),
 \end{align*}
where  $K^{(0)}: [0,1]\to \mathbb R^+ $ is a bounded function of bounded variation. The sequence $(h_n)_{n\geq 1} $ is such that $nh_n^{2d+4} \rightarrow 0$, $nh_n^d /  |\log( h_n )| \rightarrow \infty$.
\setcounter{ass:H}{\value{enumi}}
\end{enumerate}

Whereas (H\ref{ass:base_density}) and (H\ref{ass:kernel}) are rather classical in the kernel smoothing littereature (see the remarks just below Theorem \ref{main:result}), the following assumption is specific to our approach. We shall see in Section \ref{sec:subsec:examples-H3} that this is satisfied for densities with classical tails.

\begin{enumerate}[(H1)]
\setcounter{enumi}{\value{ass:H}}
\item \label{ass:Hcvlin_f_0_St} 
The function $q:\mathbb R^d \mapsto \mathbb R^+$ is bounded, integrable, and satisfies $\mathbb E[|\log(q(X_1))|]<\infty$.
There exist $\beta\in (0, 1]$ and $c>0$ such that $\int_{S_t^c} f_0(x) \, \diff x \leq c t^{\beta}$ as $t\to 0$.
For any $\gamma>0$, $b_n= \gamma (nh_n^d)^{-1/\beta}$, there exists $C>0$ such that, {as $n \to \infty$},
\begin{align*}
&\sup_{x\in S_{b_n}} \sup_{u\in [-1,1]^d}  \frac{ f_0(x+h_nu)  }{f_{0}(x) } \leq C ,\qquad 
\text{ and } \qquad h_n^{d}\lambda(S_{b_n}) \to 0.
\end{align*}
\setcounter{ass:H}{\value{enumi}}
\end{enumerate}

For the sake of clarity, the (classical) assumptions dealing with the parametric model are postponed to the appendix: (A\ref{ass:parametric_model_consistency}) and (A\ref{ass:parametric_model_assymptotic_normality}). They are taken from the monographs \cite{vandervaart:1998} and \cite{newey+m:1994}, and they mainly ensure the asymptotic normality of $\hat \theta_n$ whenever $ f_0\in \mathcal P$.

\begin{theorem}\label{main:result}
Suppose that assumptions (H\ref{ass:base_density}) to (H\ref{ass:Hcvlin_f_0_St}), and (A\ref{ass:parametric_model_consistency}) are fulfilled.
\begin{enumerate}[(i)]
\item \label{main:result_h0} When $f_0\in \mathcal P$, under (A\ref{ass:parametric_model_assymptotic_normality}) and if  $(n h_n^d) \Delta_n    \to 0$, it holds that $\hat \alpha_n \rightarrow 1$, in probability.
\item \label{main:result_h1} When $f_0\notin \mathcal P$, if $\Delta _ n \to 0$ and $( \sqrt{|\log(h_n)|/nh_n^d} +h_n^2)   |\log(\Delta_n)|^{1/\beta} \to 0$, we have that $\hat \alpha_n \to 0$, in probability.
\end{enumerate}
\end{theorem}

Appendix \ref{app:proof_main_result} is dedicated to the proof of Theorem \ref{main:result}. We did not follow the approach used in \cite{olkinSpiegelman1987}, which, we believe, is unsatisfactory because they do not consider the case when $\hat \alpha_n$ lies in the border of $[0,1]$. Actually, this is not straightforwardly remedied as the event $\hat \alpha_n = 0 $ or $\hat \alpha_n = 1$ has a non-negligible probability (as illustrated in the numerical experiments in section \ref{sec:robustness}).
The smoothness assumption stated in (H\ref{ass:base_density}) and the symmetries in the kernel function ensure a control of order $h_n^2$ of the bias $f_h(x) - f(x)$, uniformly in $x\in \mathbb R^d$ (see Lemma \ref{lemma:auxiliary_bound_kernel} stated in Appendix \ref{app:analytical_results}). Such a rate could be improved by using higher order kernels but this is not necessary here.
Assumption (H\ref{ass:kernel}), (a) and (b), are borrowed from the empirical process literature; see among others \cite{nolan+p:1987,gine+g:02,einmas0}. They permit to bound, uniformly in $x\in \mathbb R^d$, the variance term $\hat f_n(x)- f_h(x)$. The fact that the kernel has a compact support can be alleviated at the price of additional technicalities in the proof and assuming that the tails of the kernel are light enough. We did not include this analysis in the paper for reasons of clarity.

For any dimension $d\geq 1$, there exists a couple of sequences 
$(h_n,\Delta_n)_{n\geq 1}$ that fulfills the restrictions \eqref{main:result_h0}, 
\eqref{main:result_h1} of Theorem~\ref{main:result} and (H\ref{ass:kernel}). For instance, the optimal bandwidth $h_n \propto n^{-1/(d+4)}$, which minimizes the asymptotic mean integrated squared error~\citep[equation (2.5)]{wand+j:1994}, and $\Delta_n = 1/n$, is one such sequence. This means that, in practice, one can choose the bandwidth according to the various methods of the literature, see e.g.~\cite{silvermanbook1998}. 

An interesting point in Theorem \ref{main:result} is the two opposite roles played by the sequence $\Delta_n$ in (\ref{main:result_h0}) and (\ref{main:result_h1}), respectively. The consistency when $f_0\in \mathcal P$ requires $\Delta_n$ to be as small as possible whereas when $f_0\notin \mathcal P$, $\Delta_n$ must not be too close to $0$. In the proof, the case $\Delta_n = 0$ (leave-one-out) as well as $\Delta_nq(X_i) = K(0)/(nh^d)$ (OS method) need to be excluded, suggesting that these other options are not consistent under our set of assumptions.

\subsection{Distributions and bandwidth sequences satisfying~(H\ref{ass:Hcvlin_f_0_St})}
\label{sec:subsec:examples-H3}

For densities $f_0$ with unbounded supports, the verification of Assumption~(H\ref{ass:Hcvlin_f_0_St}) only depends on some tail function $g_0$ associated to the density $f_0$. The meaning of this is made precise in the following proposition.

\begin{proposition}\label{prop:g_0tof_0}
 Suppose that for any $A>0$, $\inf_{\Vert x \Vert \leq A} f_0(x)>0$ and that there exists a function $g_0$ such that $f_0(x)  / g_0(x) \to 1 $ as $\Vert x \Vert \to \infty$. Suppose that $h_n\to 0 $ obeys $nh_n^d\to \infty$. If there exist $c_2>0$ and $\beta\in (0,1]$ such that
 \begin{align*}
   \int_{ g_0(x) \leq  t } g_0(x)\, \diff x \leq c_2t^\beta ,
   \qquad \text{as }t\to 0
 \end{align*}
and if for any $\gamma>0$, $b_n = \gamma (nh_n^d)^{-1/\beta}$, there exists $A>0$, $C_2>0$ such that
 \begin{align}\label{eq:double-sup-g0}
 \sup_{ \Vert x \Vert > A,\, g_0(x)  >  b_n }  \sup_{u\in  [-1,1] } \frac{g_0(x+h_nu)}{g_0(x)}\leq C_2 ,& \qquad \text{ and } \qquad h_n^d \lambda (g_0(x) > b_n) \to 0,
 \end{align}
{as $n \to \infty$,} then (H\ref{ass:Hcvlin_f_0_St}) is valid for $f_0$ with the same value of $\beta$.
\end{proposition}

The proof of Proposition~\ref{prop:g_0tof_0} is given in Appendix~\ref{app:proof_prop}. The function $g_0$ in Proposition~\ref{prop:g_0tof_0}, not necessarily a proper density function,  represents the rate of decrease of $f_0(x)$ as $\Vert x \Vert \to \infty$. 

\begin{example}[Mixture of densities]\label{ex:mixture-of-densities}
  Let $d=1$. Let $f_0(x) = \pi_1 f_1(x) + \pi_2 f_2(x)$, $\pi_1 > 0$, $\pi_2 > 0$, 
  $\pi_1 + \pi_2 = 1$, where $f_1$ and $f_2$ are densities such that 
  $f_1(x) / f_2(x) \to 0$ as $|x| \to \infty$. Take $g_0(x) =  \pi_2 f_2(x)$. Then,  as $|x| \to \infty$,
  \begin{align*}
    \frac{ f_0(x) }{ g_0(x) } =
     \frac{\pi_1 f_1(x) }{  \pi_2 f_2(x) } + 1 \to 1.
  \end{align*}
  Hence the verification of~(H\ref{ass:Hcvlin_f_0_St}) by $f_0$ only depends on the component $f_2$. 
\end{example}

Putting $g_0 \propto f_0$ (the symbol $\propto$ stands for proportionality) in Proposition~\ref{prop:g_0tof_0} amounts to check~(H\ref{ass:Hcvlin_f_0_St}) directly, which is done in the following examples.

\begin{example}[Gaussian tails]
Let $d=1$ and $g_0(x) =\kappa_1 \exp(-\kappa_2x^2)$, with $\kappa_1>0$, $\kappa_2>0$. For clarity the computations are provided for $\kappa_1 = \kappa_2 =1$ but can easily be extended for arbitrary values. Because $  \int_{\exp(-x^2) \leq t } \exp(-x^2) \, \diff x \leq t $, as $t\to 0$, we have that $\beta= 1$. Moreover, for $0<b_n<1$, we have
\begin{align*}
  \sup_{\exp(-x^2)  >  b_n }  \sup_{u\in  [-1,1] } \exp(-(x+h_nu)^2 + x^2) 
  &\leq   \sup_{\exp(-x^2)  >  b_n }  \sup_{u\in  [-1,1] } \exp(-2h_nxu)\\ 
  & = \sup_{\exp(-x^2)  >  b_n } \exp(2h_n|x| )\\
  &\leq    \exp(2 \sqrt {- \log(b_n) } h_n )
\end{align*}
Therefore, a sufficient condition on $h_n$ guaranteeing~\eqref{eq:double-sup-g0}
is that $h_n^2 \log(n)\to 0$, which is satisfied under~(H\ref{ass:kernel}).
\end{example}

\begin{example}[Exponential tails] Let $d=1$ and $g_0(x) = \kappa_1 \exp(-\kappa_2 x)$, with $\kappa_1>0$, $\kappa_2>0$.
The computations are very similar to the one presented in the Gaussian case. We find $\beta = 1$ and the condition on $h_n$ becomes $h_n\log(n) \to 0$ which is always true under (H\ref{ass:kernel}). Hence, as for Gaussian tails, when the tails are exponential, (H\ref{ass:Hcvlin_f_0_St}) is automatically satisfied under (H\ref{ass:kernel}).
\end{example}

\begin{example}[Polynomial tails]
Let $d=1$ and $g_0(x) = \kappa_1|x|^{-k}$ with $\kappa_1>0$, $k>1$. For simpicity, as in the Gaussian example, we focus on $\kappa_1 = 1$. We find that $\beta = (k-1)/k$. For $h_n<|A|$, we have
\begin{align*}
\sup_{|x|> A,\, |x | \leq b_n^{-1/k} } \sup_{u\in[-1,1]^d} \frac{|x|^k}{|x+h_nu|^k} = \sup_{|x|> A,\, |x | \leq b_n^{-1/k} } \frac{|x|^k}{(|x|-h_n)^k} =
\frac{1}{(1-h_n/A)^k} \overset{n\to \infty}{\to} 1.
\end{align*}
Finally, since 
$h_n \lambda( g_0>b_n ) =  2 b_n^{-1/k} = 2 \gamma^{-1/k} (nh_n)^{1/(\beta k)}$, a sufficient condition on $h_n$ guaranteeing~\eqref{eq:double-sup-g0}
is that $nh_n^k \to 0 $. 
\end{example}

The three examples considered above are informative on the interplay between the tails of $f_0$ and the choice of $h_n$. For distribution with light enough tails, including Gaussian, exponential and polynomials with $k\geq 6$, the conditions on $h_n$ required by (H\ref{ass:Hcvlin_f_0_St}) are already fulfilled when assuming (H\ref{ass:kernel}). Consequently, the optimal bandwidth which has order $n^{-1/5}$ is included by our set of assumptions. In contrast, as soon as $k<6$ in the polynomial case, we have the additional condition that $nh_n^k \to 0$. 


\section{Numerical illustrations}\label{sec:numerical_results}


In all the simulation experiments, we have set 
$\Delta_n=1/n$, $K(x) = \exp(-x^2/2)/\sqrt{2\pi}$ and $q(x)=t_\nu((x-\mu_q)/\sigma_q)$ where $t_\nu$ is the density of a Student-t distribution with $\nu=3$ degrees of freedom, $\mu_q=0$ and $\sigma_q=100$. With such a large variance and heavy tails, this choice of $q$ is non informative. We made $\mu_q$ and $\sigma_q$ vary but the results were very similar, suggesting that the choice for $q$ has little effect in practice (at least for light-tailed distribution). In all the experiments but those in Section~\ref{sec:robustness}, the bandwidth was chosen according to the well known rule of thumb given in~\cite{silvermanbook1998} (p. 48, equation (3.31)). The choice of the bandwidth is discussed in Section~\ref{sec:robustness}. In Section~\ref{sec:experiment}, we study the behavior of the fitness coefficient and the performance
of the estimators with respect to the amount of evidence of the model.
In Section~\ref{sec:example in bivariate density estimation}, we use the LR method for protection against misspecification.
All the numerical experiments were carried out with the \verb+R+ software.

\subsection{Sensitivity to the bandwidth: comparison of the fitness coefficient and the OS coefficient}\label{sec:robustness}

In this section, we study how a change in the bandwidth affects the fitness coefficient and the OS coefficient. 
We reanalyze the data of Olkin and Spiegelman~\cite{olkinSpiegelman1987},  consisting of yearly wind speed maxima taken in the north direction in Sheridan, Wyoming. There are 20 observations for the years 1958 to 1977: 70, 61, 61, 60, 61, 63, 61, 67, 61, 62, 47, 67, 61, 49, 55, 65, 57, 51, 47, 56. 
The parametric model is a Gumbel model, that is, $\log f_\theta(x) = (x-\mu)/\sigma - \exp((x-\mu)/\sigma)$, where $\theta=(\mu,\sigma)$, $\mu$ is a real location parameter and $\sigma>0$ a dispersion parameter obeying $\var_{f_\theta}=\pi^2\sigma^2/6$ and
$\E_{f_\theta}=\mu-\sigma\gamma$, where $\gamma \approx 0.58$ is the Euler-Mascheroni constant. The maximum likelihood estimator is given by $\hat{\theta}_n \approx (62.1, 5.4)$.

\begin{figure}[ph]
  \centering
  \subfloat[]{
    \label{fig:subfig:robustness}
    \includegraphics[width=.5\textwidth]
    {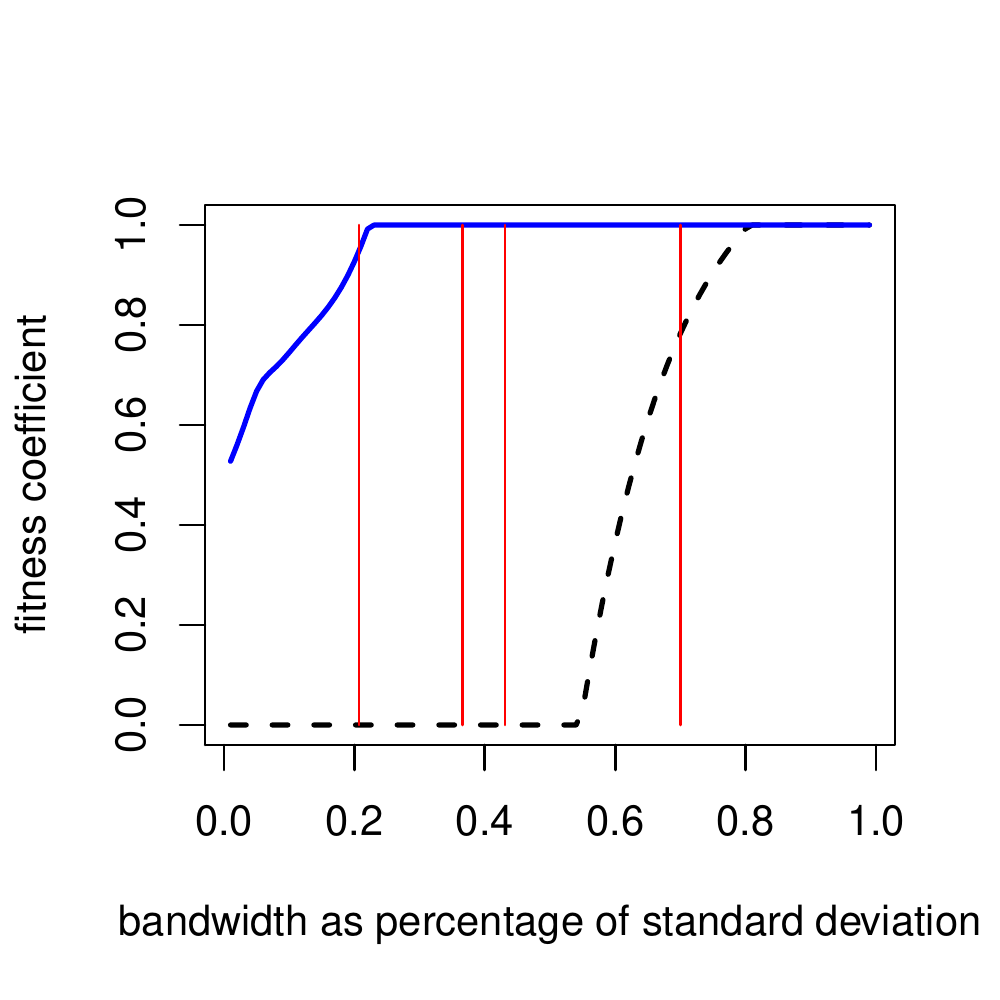} }
  \subfloat[]{
    \label{fig:subfig:touchup}
    \includegraphics[width=.5\textwidth]
    {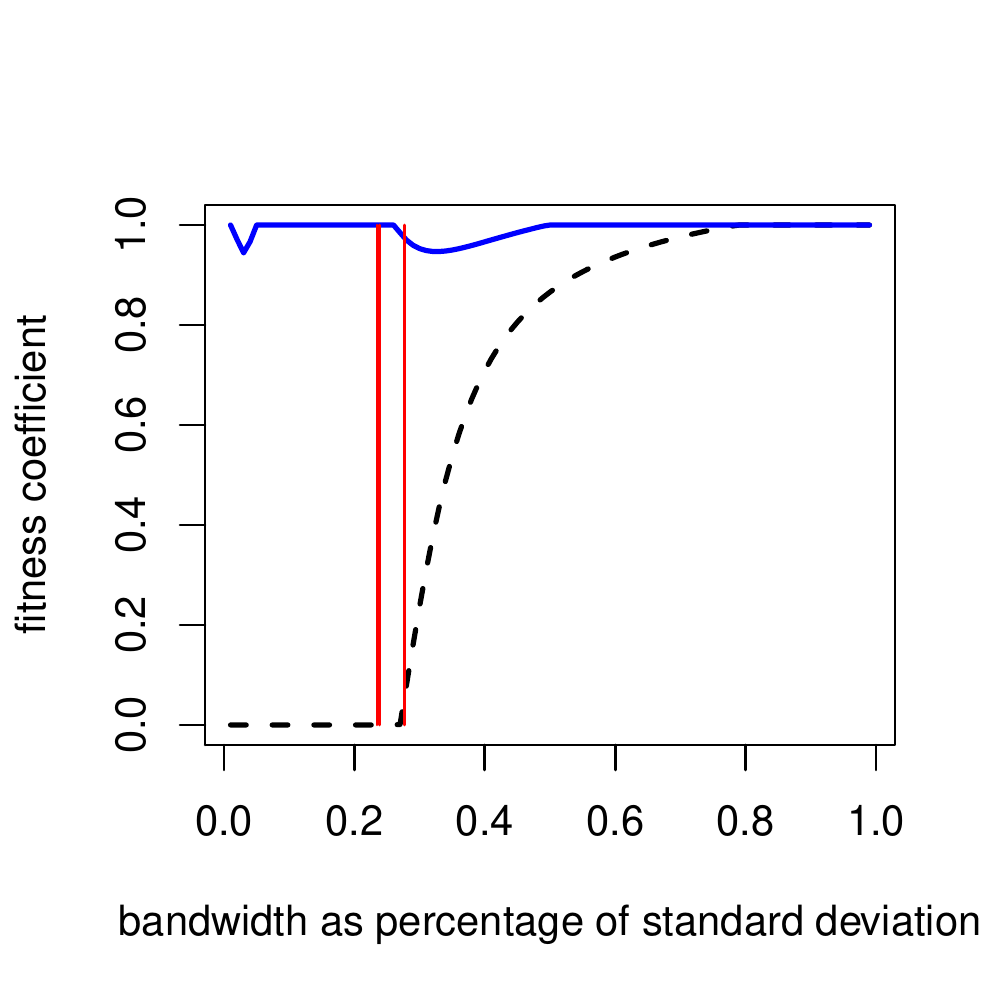} }\\    
  \subfloat[]{
    \label{fig:subfig:touchdown}
    \includegraphics[width=.5\textwidth]
    {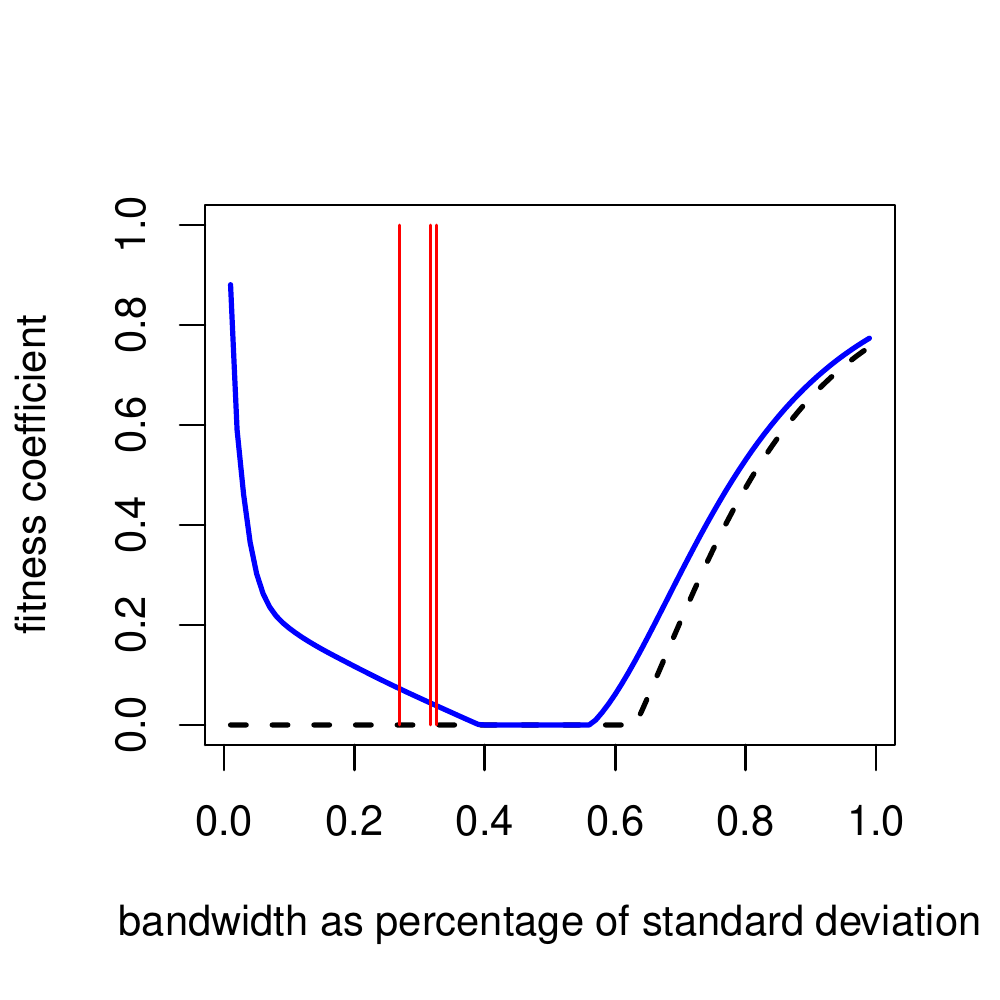} }
  \subfloat{
    \includegraphics[width=.5\textwidth]
    {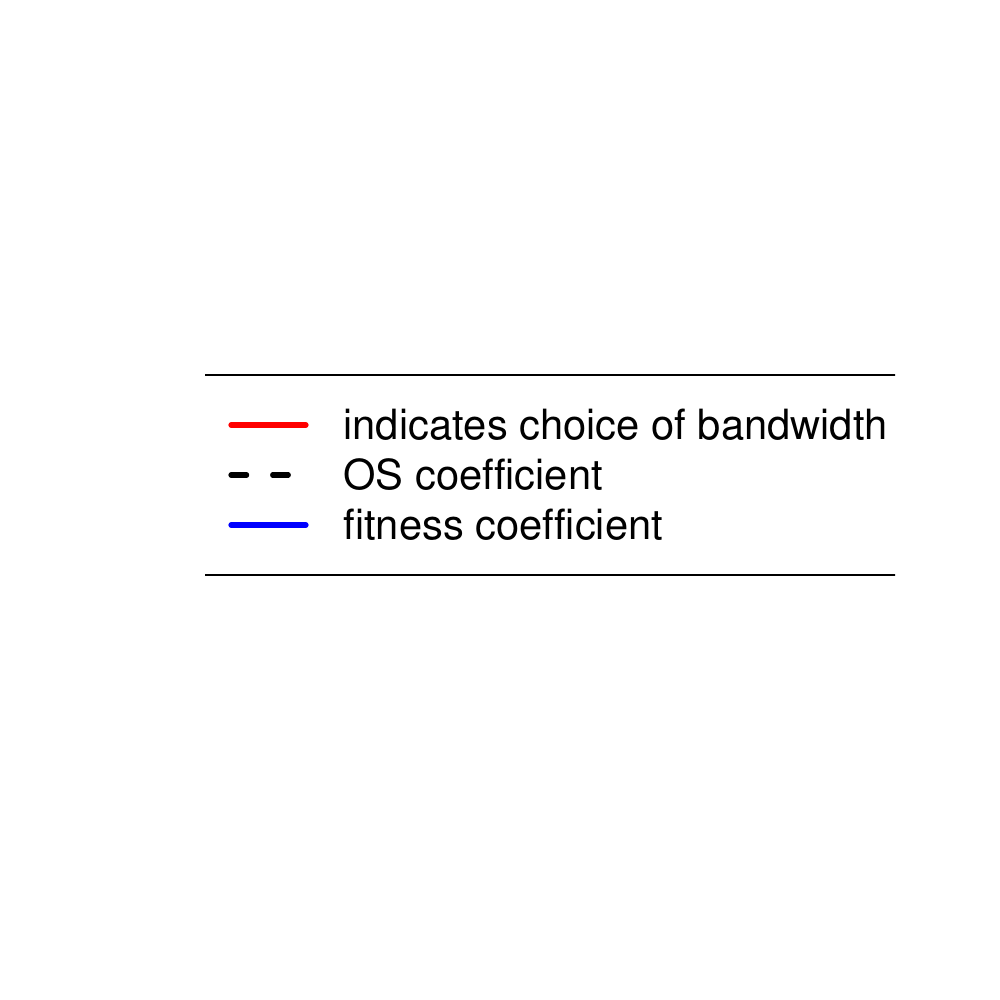} }
  \caption{Values of the fitness coefficient and the OS coefficient as a function of the bandwidth $h$, expressed as a proportion of the standard error of the data. Plain blue line: fitness coefficient. Dashed black line: OS coefficient. The red sticks indicate various  bandwidth values chosen according to the literature
 (see text); \protect\subref{fig:subfig:robustness} wind speed data, 
 \protect\subref{fig:subfig:touchup} simulations under a Gumbel model, 
\protect\subref{fig:subfig:touchdown} simulations under a Gaussian model. In all cases the fitted model is a Gumbel model.}
  \label{fig:alpha-as-function-of-bandwidth}
\end{figure}

Let $h$ denote the bandwidth. In~\cite{olkinSpiegelman1987}, it was arbitrarily chosen $h=0.7s$, where $s$ is the standard deviation of the data. This yields $\alphahat^{\text{OS}} \approx 0.8$. But if $h \approx 0.43s, \, h \approx 0.37s$ or 
$h \approx 0.21s$ then one gets $\alphahat^{\text{OS}} \approx 0$. All the above values for $h$ are grounded by well-known bandwidth selection methods, see the textbook~\cite{silvermanbook1998} (p. 47, eqn (3.30) and p. 48, eqn (3.31)) 
and~\cite{sheatherJones1991bandwidth}. By contrast, the fitness coefficient yields $\alphahat \approx 1$. These findings are summarized in Figure~\ref{fig:alpha-as-function-of-bandwidth}~\subref{fig:subfig:robustness}
, where the coefficients are represented as functions of $h$. We see that the OS coefficient is sensitive to the choice of the bandwidth: a slight difference in $h$ can yield a large difference in $\alphahat^{\text{OS}}=\alphahat^{\text{OS}}(h)$ especially in the range $0.4\le h\le 0.8$. On the opposite, the fitness coefficient is more robust: the estimated value for $\alphahat(h)$ remains close to one in a large range for $h$. In Figure~\ref{fig:alpha-as-function-of-bandwidth}~\subref{fig:subfig:robustness}, the fitness coefficient and the OS coefficient contradict each other and no more credit can be given to any one of them because the ground truth is unknown.


To observe the behavior of the coefficients when the model is known to be true,
we simulated $n=400$ observations according to a Gumbel distribution with mean  and standard deviation equal to those of the wind speed data, that is, 59.1 and 6.55 respectively.
The results are shown in 
Figure~\ref{fig:alpha-as-function-of-bandwidth}~\subref{fig:subfig:touchup}.
One has $\alphahat(h)\approx 1$ whatever $h$ while 
 $\alphahat^{\text{OS}}(h)\le 0.2$ for all $h$ chosen by 
the bandwidth selection methods of the literature. These results tend to indicate that the fitness coefficient is consistent but the OS coefficient is not.
Let us note that Figure~\ref{fig:alpha-as-function-of-bandwidth}~\subref{fig:subfig:robustness} and~\subref{fig:subfig:touchup} are similar, making the Gumbel model  plausible. The difference spotted in the range $0 \le h \le 0.2$ can be explained by the ties of the wind speed data. (When $h$ is small, one can see that~(\ref{eq:our method}) is close to the likelihood of a Bernouilli trials experiment, the maximizer of which is given by the proportion of untied observations, here one half.) 

Whenever the model is wrong, we found on simulations that for most reasonable (that is, found in the literature as above) values of $h$, the values of the coefficients are close to  zero, as expected. This is illustrated in Figure~\ref{fig:alpha-as-function-of-bandwidth}~\subref{fig:subfig:touchdown}:
the model is still Gumbel, but the $n=400$ data points were generated according to a Gaussian distribution with mean 59.1 and standard deviation 6.55. 

\subsection{Performance of the methods when the model and the truth intertwine}\label{sec:experiment}
Parametric estimators perform better than kernel density estimators when the model is approximately true, but worse otherwise. Can the semiparametric combination be uniformly best? Does the fitness coefficient goes to unity as the model approaches the truth? 

To get some insight, the following numerical experiment was done. We generated samples of size $n=400$ according to a density $f_t$, for several values $t$ in a certain index set, representing the ``distance'' between $f_t$ and the model. Two settings have been tested.

\begin{description}

\item[Setting 1] The parametric model is given by 
 $f_\theta \sim N(\theta, 1)$ and the curve of true distributions is given by $f_t \sim N(0, (1+t)^2)$. The intersection between the model $\{f_\theta\}$ and the family $\{f_t\}$ is given by
$f_0 \sim N(0,1)$; that is, $\theta=t=0$.

\item[Setting 2] The parametric model is given by
 $f_\theta \sim N(0, \theta^2)$ and the curve of true distributions is given by $f_t \sim N(t, 1)$. The intersection between the model $\{f_\theta\}$ the family $\{f_t\}$ is given by
a $N(0,1)$ as well.

\end{description}

For each $t$, we compute the maximum likelihood estimator, the standard kernel density estimator, the fitness coefficient, the OS coefficient, and the semiparametric density estimator. The semiparametric density estimator is the combination between the maximum likelihood estimator and the kernel density estimator where the mixing coefficient can be either the fitness coefficient (LR method) or the OS coefficient (OS method). To assess the performances of the estimators, we compute the L2-distance to $f_t$. The above procedure is repeated 500 times so that the errors are averaged over the repetitions. 

\begin{figure}[phb]
 \subfloat[]{
   \label{fig:subfig:errors-n400nrep500setting1seed1}
   \includegraphics[width=.5\textwidth]
   {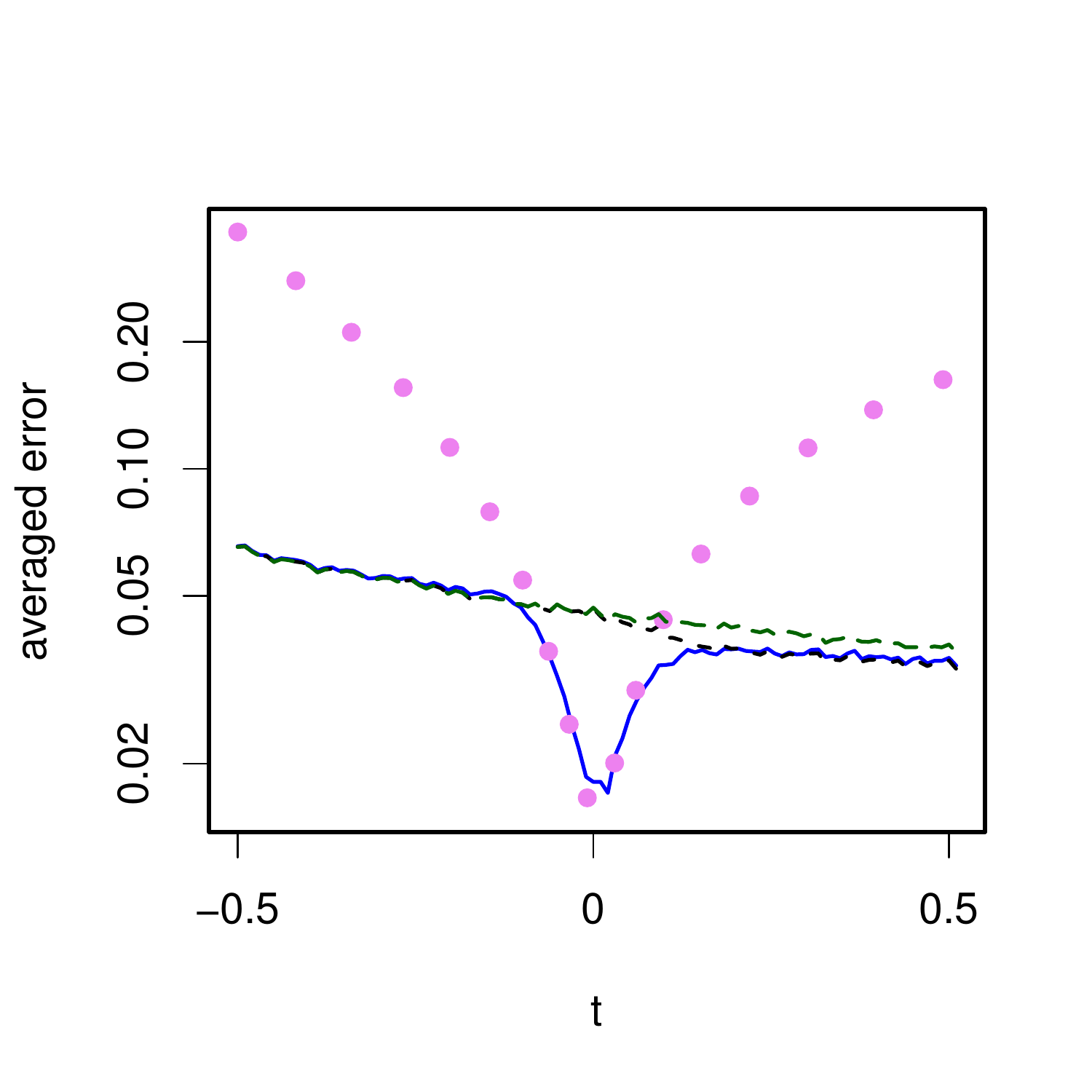} }
 \subfloat[]{
   \label{fig:subfig:mixingCoefficients-n400nrep500setting1seed1}
   \includegraphics[width=.5\textwidth]
   {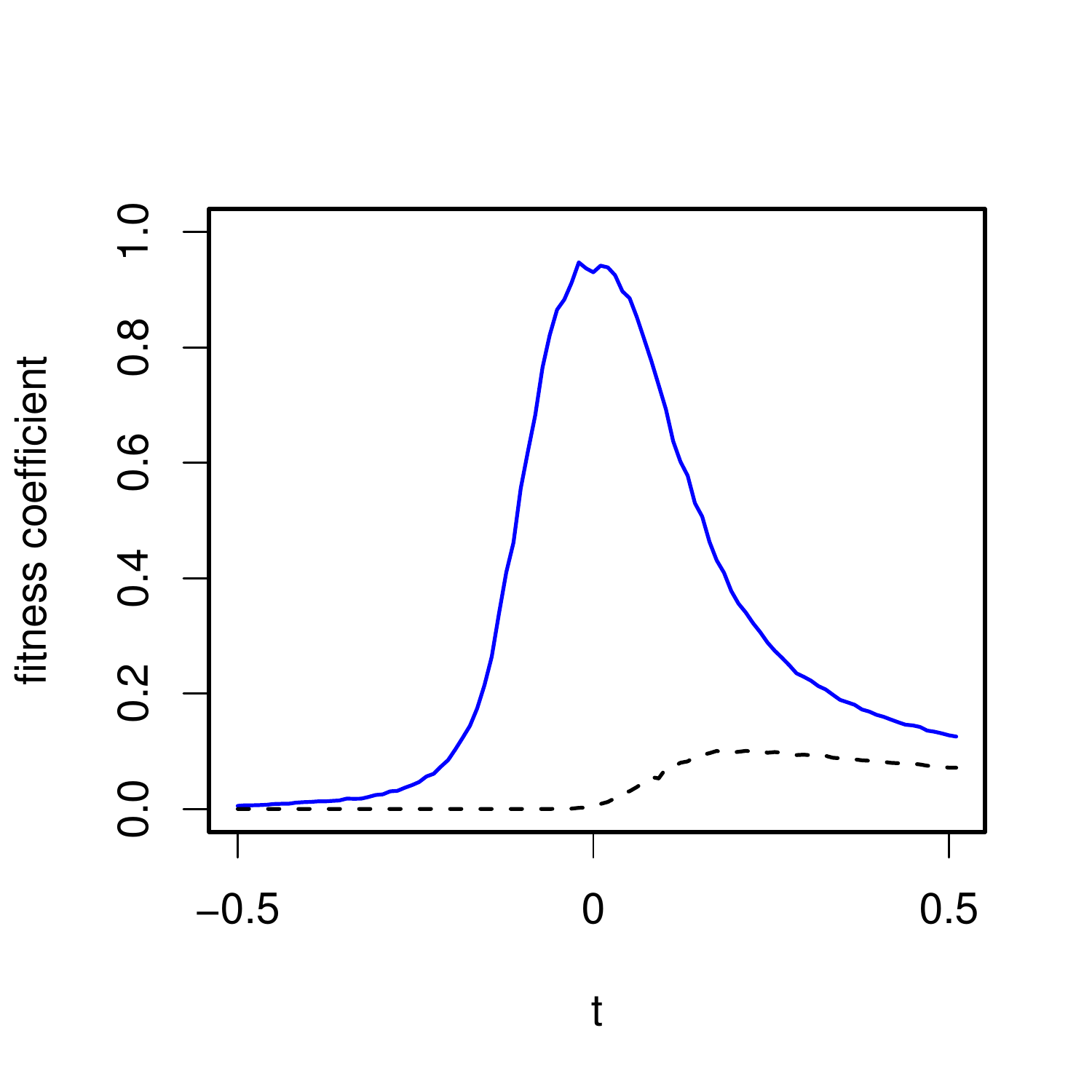} }\\
 \subfloat[]{
   \label{fig:subfig:integratedErrors-n400nrep500setting1seed1}
   \includegraphics[width=.5\textwidth]
   {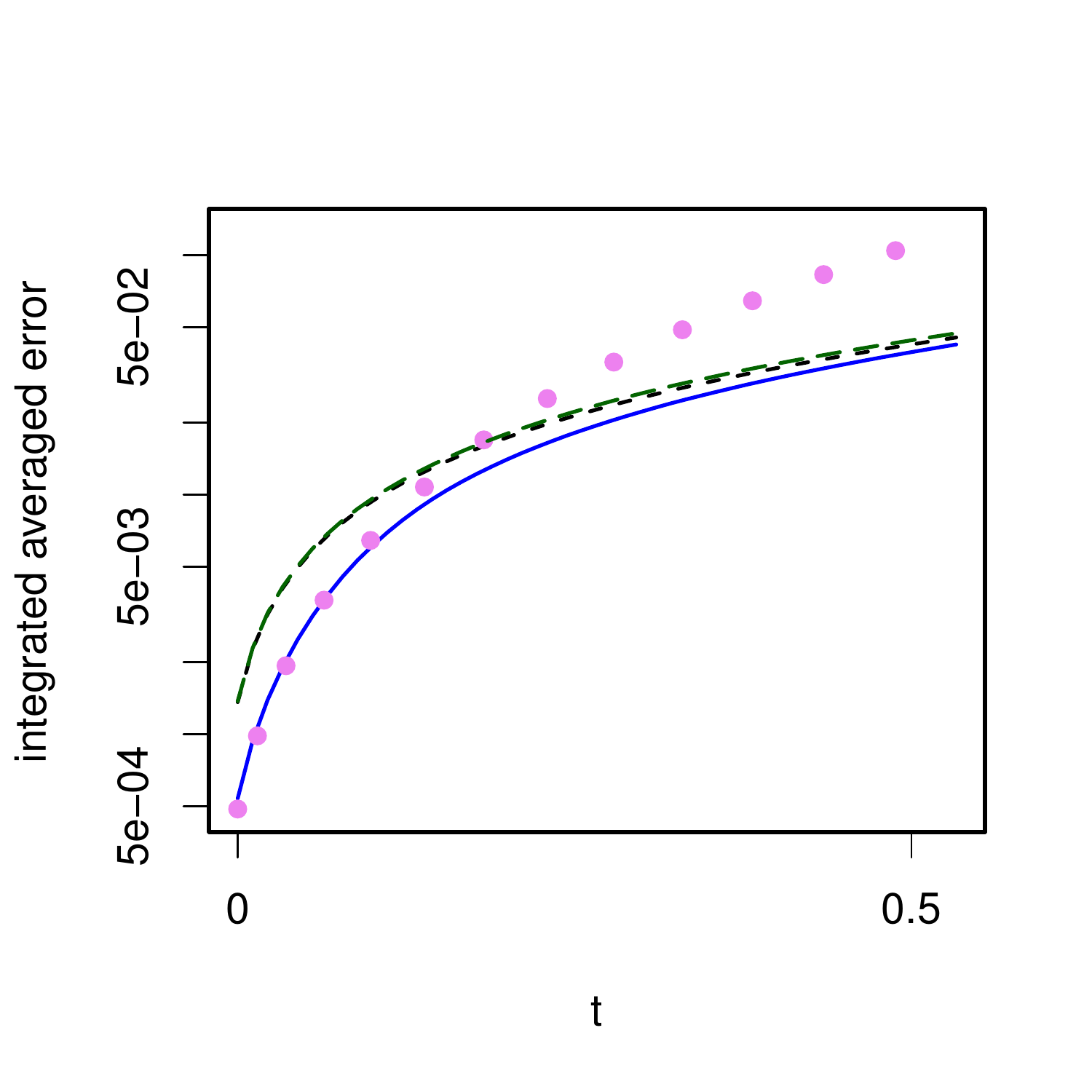} }
 \subfloat{
   \includegraphics[width=.5\textwidth]
   {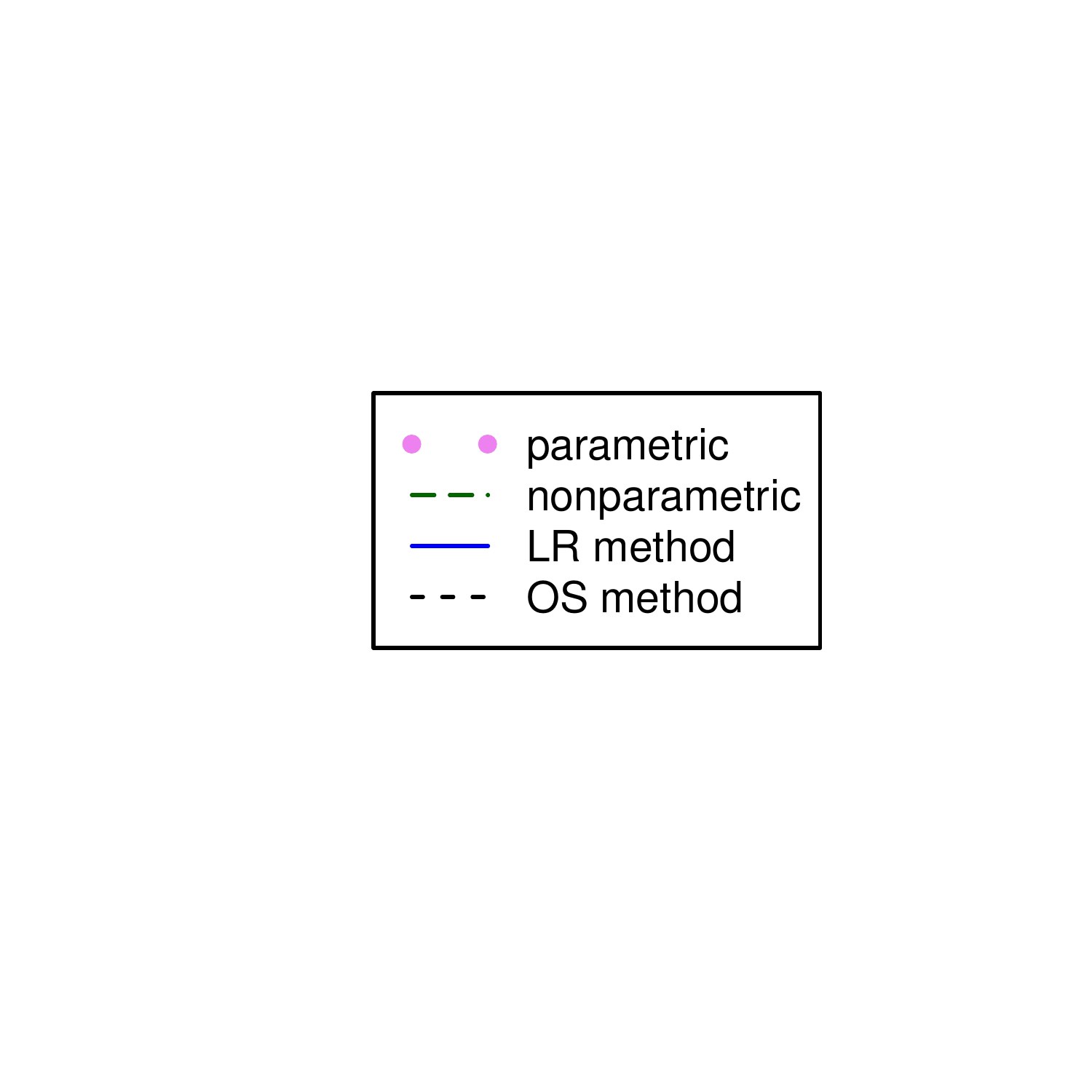} }
 \caption{Performance of the methods when the truth 
   $\{N(0,(1+t)^2), \, -0.5<t<0.5\}$ approaches the model 
   $\{N(\theta,1), \, -\infty<\theta<\infty\}$
   until they intersect at $t=0$. The L2-distance averaged over the replications is pictured
   in~\protect\subref{fig:subfig:errors-n400nrep500setting1seed1} for the parametric estimator, the nonparametric estimator, the OS method and the LR method.
   The integrated averaged distance is pictured 
   in~\protect\subref{fig:subfig:integratedErrors-n400nrep500setting1seed1}. Figure~\protect\subref{fig:subfig:mixingCoefficients-n400nrep500setting1seed1} pictures the values of the fitness coefficient and the OS coefficient averaged over the replications.}
 \label{fig:n400setting1}
\end{figure}

Figure~\ref{fig:n400setting1} summarizes the results for the first setting. 
The errors for the parametric estimator, shown in Figure~\ref{fig:n400setting1}~\protect\subref{fig:subfig:errors-n400nrep500setting1seed1}, shrink sharply as the model and the truth intersect. The error for the nonparametric estimator is approximately constant. We see that the OS method performs poorly: it fails to give accurate estimates near the truth. This behavior is explained in Figure~\ref{fig:n400setting1}~\protect\subref{fig:subfig:mixingCoefficients-n400nrep500setting1seed1}, where we see that the values of the OS coefficient barely exceed 0.1. This is not the case for the fitness coefficient; the values stretch entirely the range $[0,1]$ and {are consistent} with the proximity between the truth and the parametric model.
 As a consequence, coming back to Figure~\ref{fig:n400setting1}~\protect\subref{fig:subfig:errors-n400nrep500setting1seed1}, the error of the LR method is near the minimum of
the parametric and nonparametric errors. This means that, in practice, however close our parametric model is to the truth, we never lose by choosing the LR method. Even more interestingly is the fact that in the region where the parametric and the nonparametric estimators perform similarly, the LR method performs better: this corresponds to the values $t\approx -0.10$ and $t\approx 0.15$. This fact is clearly seen in Figure~\ref{fig:n400setting1}~\protect\subref{fig:subfig:integratedErrors-n400nrep500setting1seed1} which pictures the averaged error integrated in the interval $[-t,t]$: the LR method always has the lowest curve.

The results for $n=50, 100, 200$ and for setting 2 are similar and not shown here to limit the length of the paper.

\subsection{Application to multivariate density estimation}\label{sec:example in bivariate density estimation}
It is well known that building accurate multivariate parametric models is an uncertain and difficult task. One way of addressing this problem consists of decomposing the target density $f_0$ into a copula $c$ and the marginal densities
$f_1,\dots,f_d$, that is,
\begin{align*}
 f(x_1, \dots, x_d)
 = c( F_1(x_1), \dots, F_d(x_d) ) f_1(x_1) \cdots f_d(x_d)
\end{align*}
(here the $\{F_j\}$ stand for the distribution functions). This decomposition, also known as Sklar's theorem, is unique provided that the $\{F_j\}$ are continuous; for more details about copulas, see e.g.~\cite{genest2007everything} or the books~\cite{nelsen_introduction_2006, joeBook2014}.
The copula is assumed to belong to a parametric model 
$\{c_\xi,\,\xi\in\Xi\}$ and the true underlying parameter $\xi$ is estimated~\cite{genest1995semiparametric} by 
\begin{align*}
 \hat{\xi}
 = \underset{\xi\in\Xi}{\arg\max}\, 
 \sum_{i=1}^n \log c_{\xi}\left( \frac{R_{i,1}}{n}, \dots, \frac{R_{i,d}}{n} \right),
\end{align*}
where $R_{i,j}$ is the rank of $X_{i,j}$ among $(X_{1,j},\dots,X_{n,j})$ and $X_{i,j}$ stands for the $j$-th coordinate of the $i$-th observation. The marginals are estimated in a separate step. If one of the marginals is misspecified, the estimation of the joint distribution is biased. In the following, a computer experiment illustrates that the LR method can help to 
reduce this bias by avoiding misspecification.


\begin{figure}
\centering
 \subfloat[]{
   \label{fig:subfig:n200-script4-marginal1}
   \includegraphics[width=.4\textwidth]
   {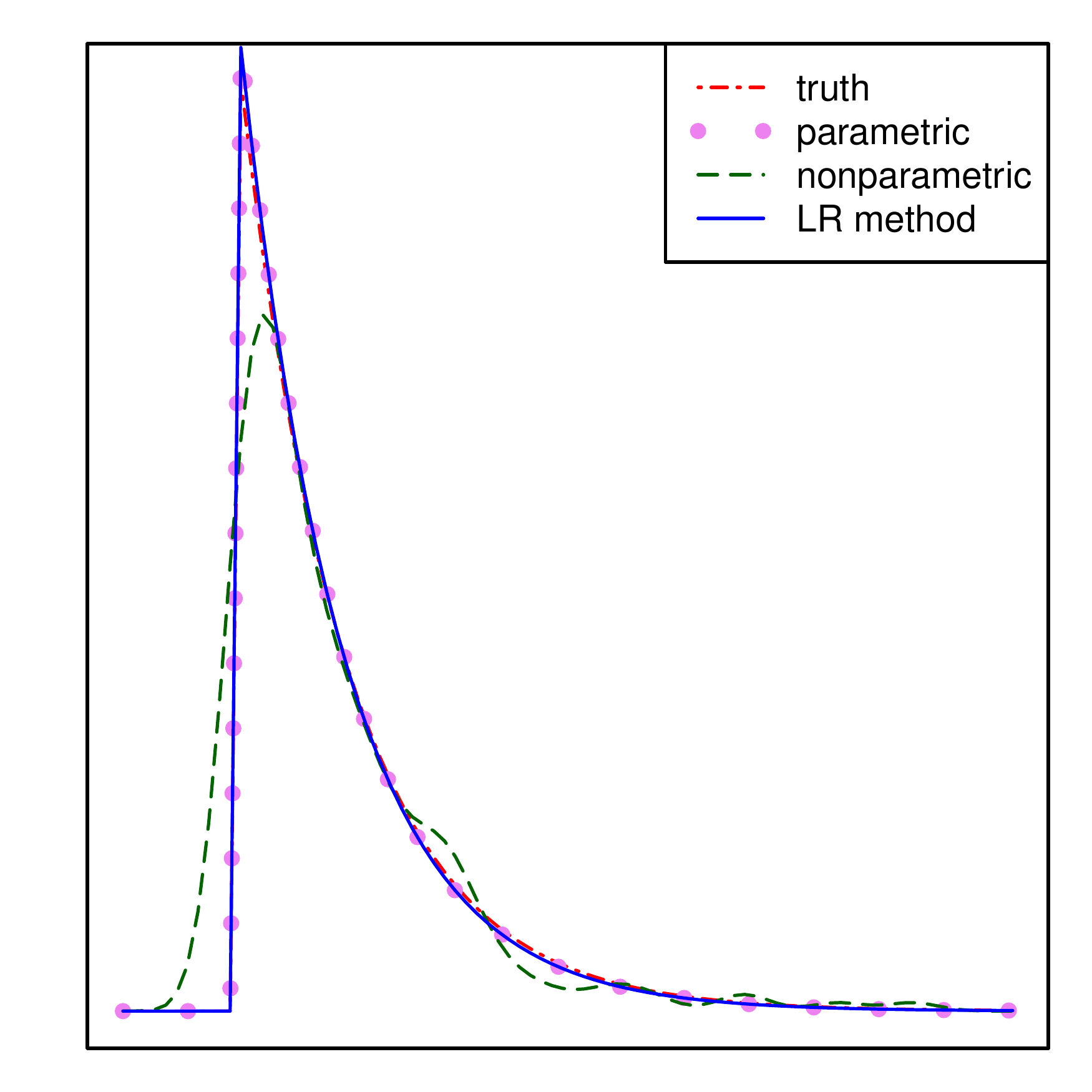} }
 \subfloat[]{
   \label{fig:subfig:n200-script4-marginal2}
   \includegraphics[width=.4\textwidth]
   {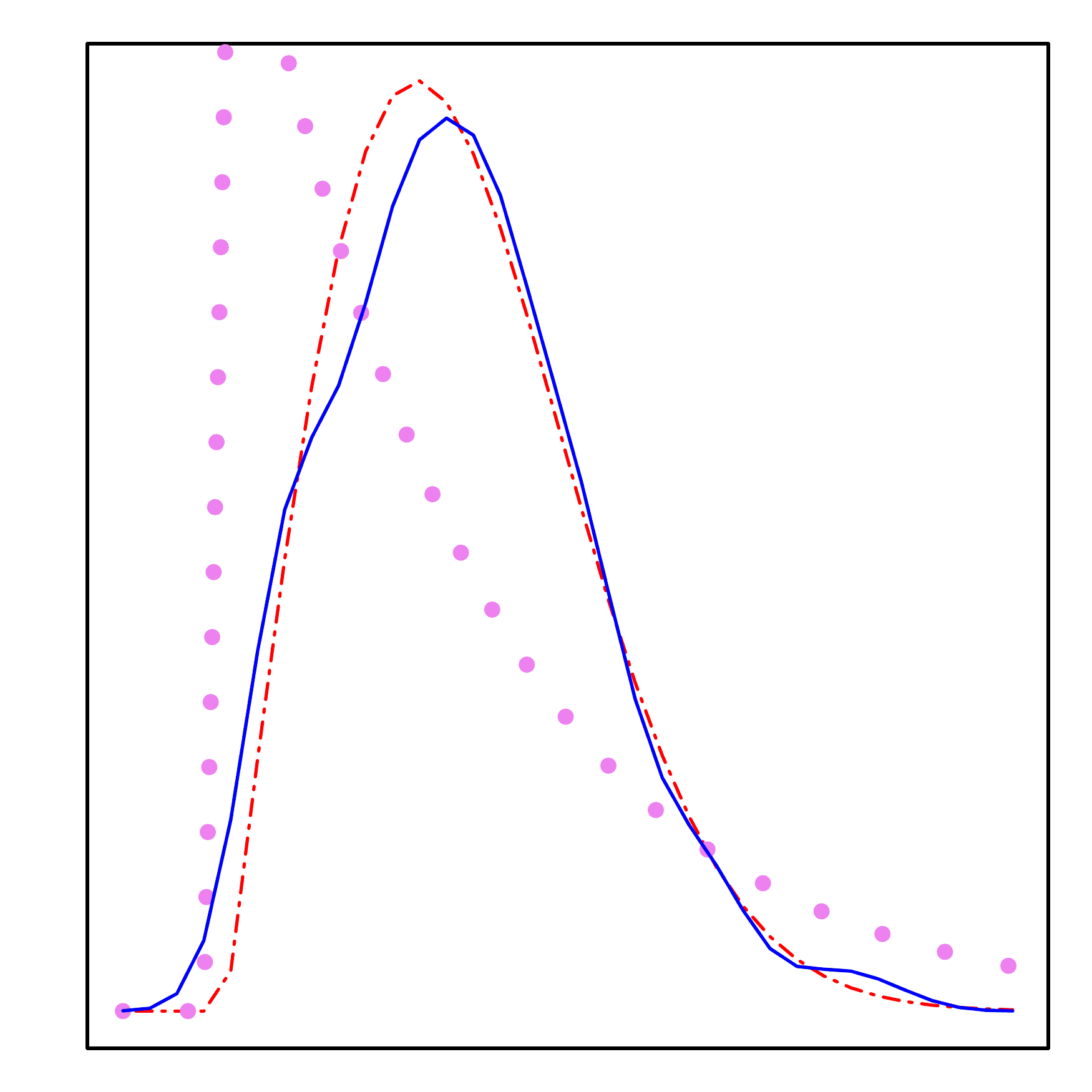} }
 \caption{Estimated densities for
   \protect\subref{fig:subfig:n200-script4-marginal1}
   the first marginal and
   \protect\subref{fig:subfig:n200-script4-marginal2}
   the second. The dotted violet, dashed green and plain blue lines
 are the parametric, nonparametric and the
 semiparametric estimates, respectively.
The dotted-dashed red line corresponds to the
true density. The size of the dataset is $n=200$.}
 \label{fig:n200-script4-marginals}
\end{figure}

We have generated datasets of size $n=25, 50, 100, 150, \dots, 500$ with a copula of the form (a so-called Gumbel copula)
\begin{align}
 \label{eq:gumbel copula}
 C_\xi( u_1, u_2 ) 
 = \exp\left\{ - \left[ (-\log u_1)^\xi + (-\log u_2)^\xi \right]^{1/\xi} \right\},
 \qquad \xi \ge 1,
\end{align}
with $\xi=3$ and marginals 
$f_1 \sim E(2)$, $f_2 \sim W(2,{1}/{2})$ where $E(\lambda)$ is an exponential distribution with mean ${1}/{\lambda}$ and $W(a,b)$ is a Weibull distribution with shape $a>0$ and scale $b>0$, that is,
\begin{align*}
 f_2(x; a,b)
 = \frac{a}{b} \left(\frac{x}{b}\right)^{a-1} \exp\left( -\left(\frac{x}{b}\right)^a \right), 
 \qquad x>0.
\end{align*}
For each of the simulated datasets, the copula parameter $\xi$ was estimated as mentioned above and the marginals were estimated under three scenarios. In the first scenario we estimate them nonparametrically with the standard kernel density estimator. In the second scenario we do as if both marginals were exponentially distributed and compute the maximum likelihood estimator. In the third scenario we form the convex combination with the maximum likelihood estimator and the standard kernel density estimators, where the mixing coefficient is the fitness coefficient (LR method).

\begin{figure}
  \centering
  \subfloat[Truth]{
    \label{fig:subfig:n200-script4-truth}
    \includegraphics[width=.35\textwidth]
    {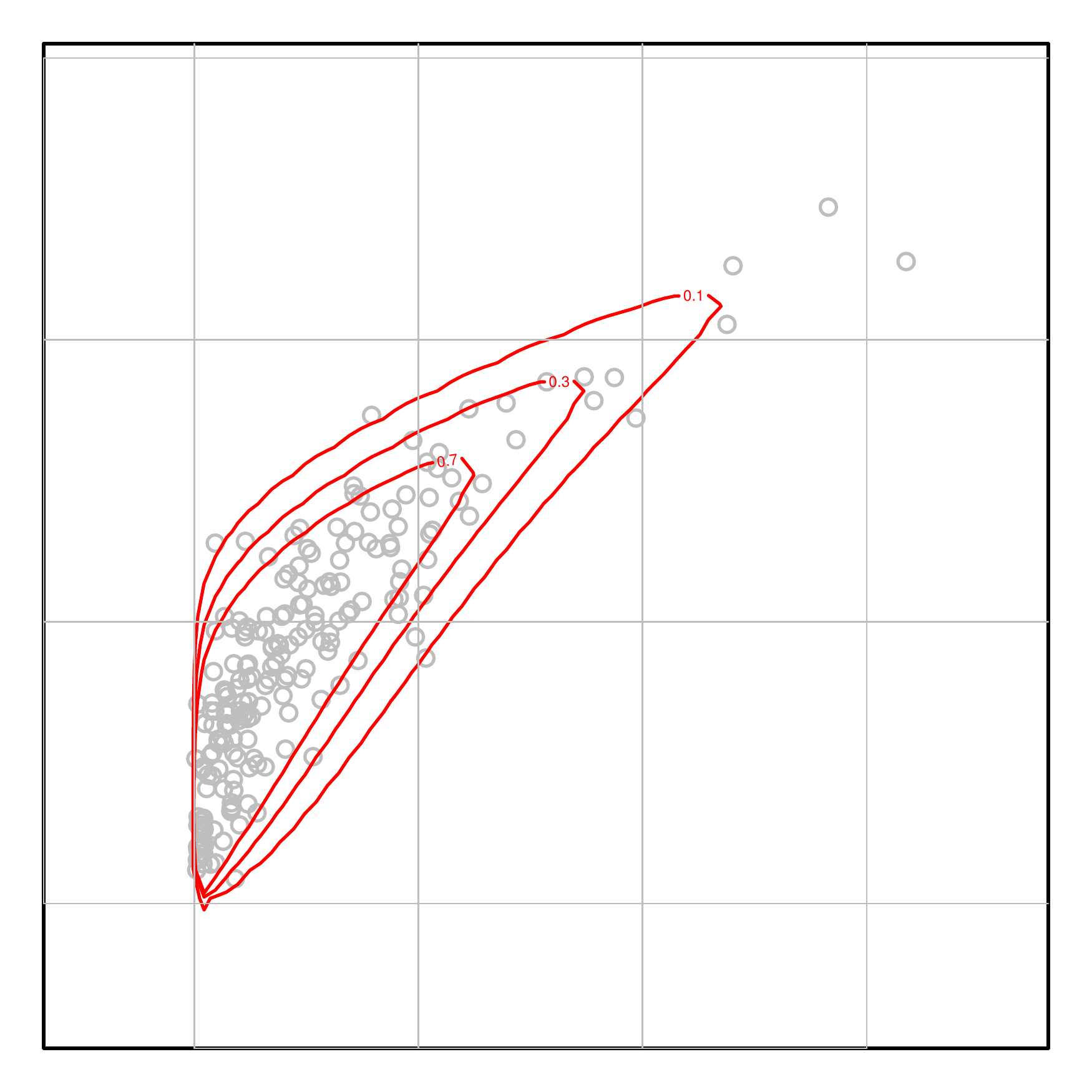} }
  \subfloat[Parametric marginals]{
    \label{fig:subfig:n200-script4-para}
    \includegraphics[width=.35\textwidth]
    {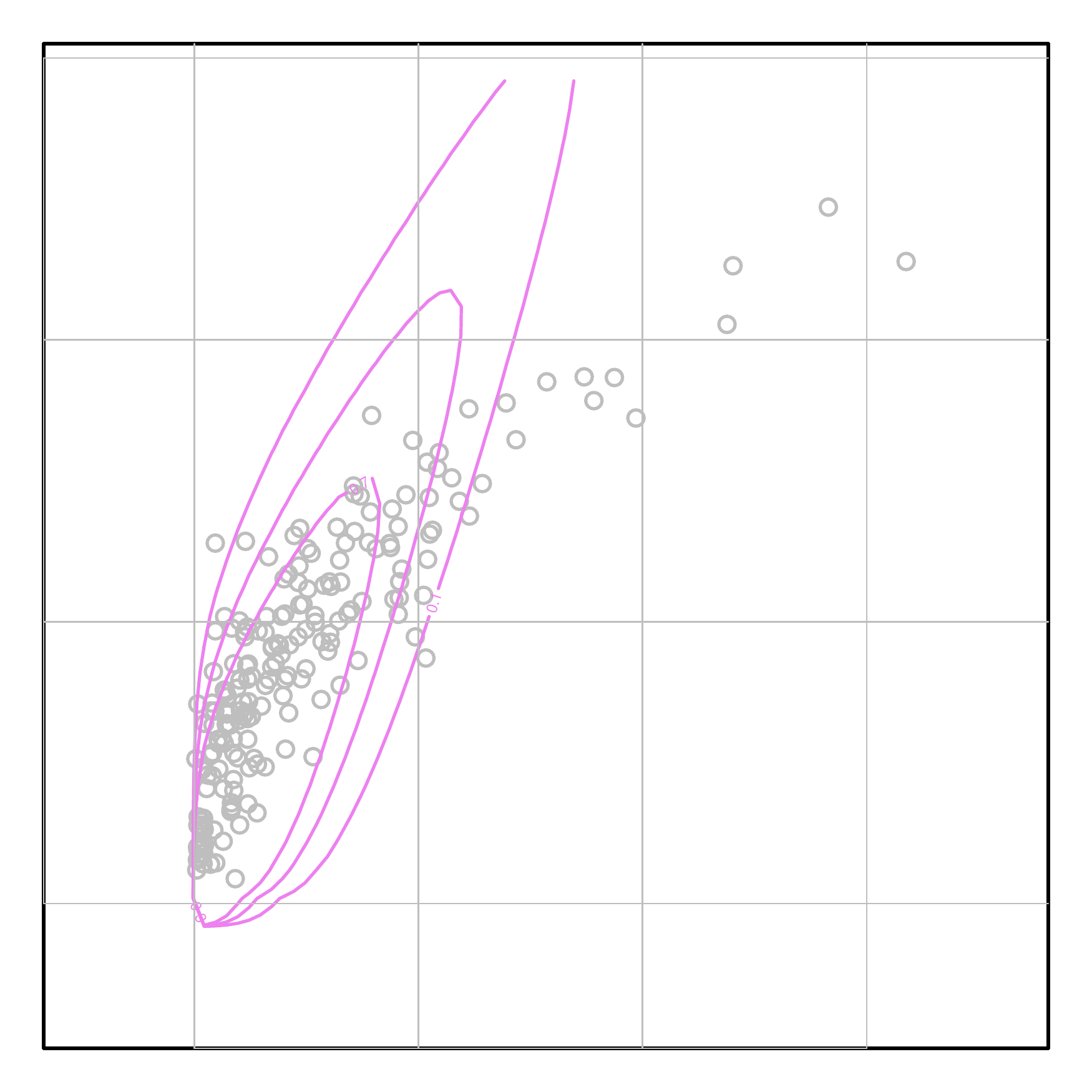} }\\
  \subfloat[Nonparametric marginals]{
    \label{fig:subfig:n200-script4-nonpara}
    \includegraphics[width=.35\textwidth]
    {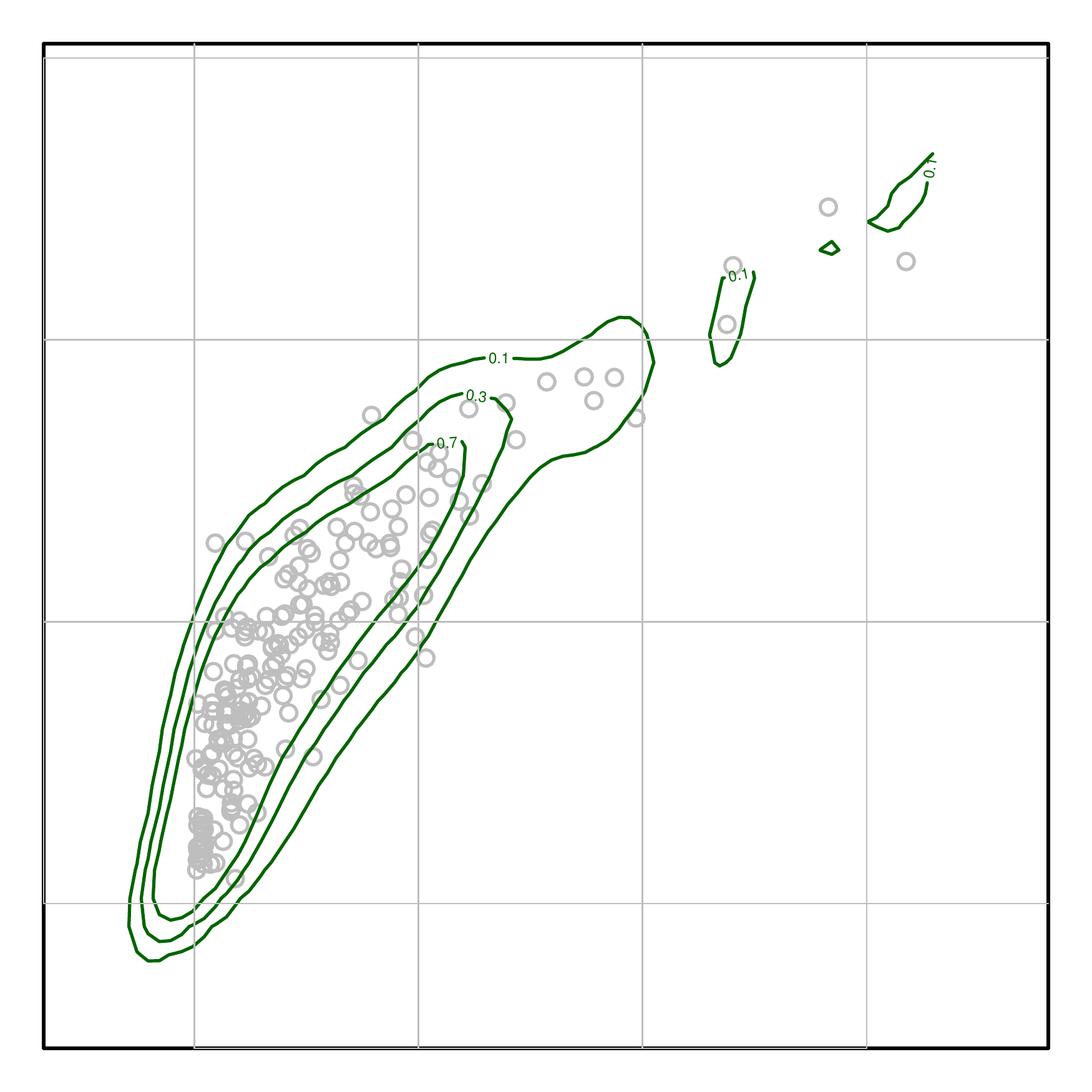} }
  \subfloat[Semiparametric marginals]{
    \label{fig:subfig:n200-script4-semi}
    \includegraphics[width=.35\textwidth]
    {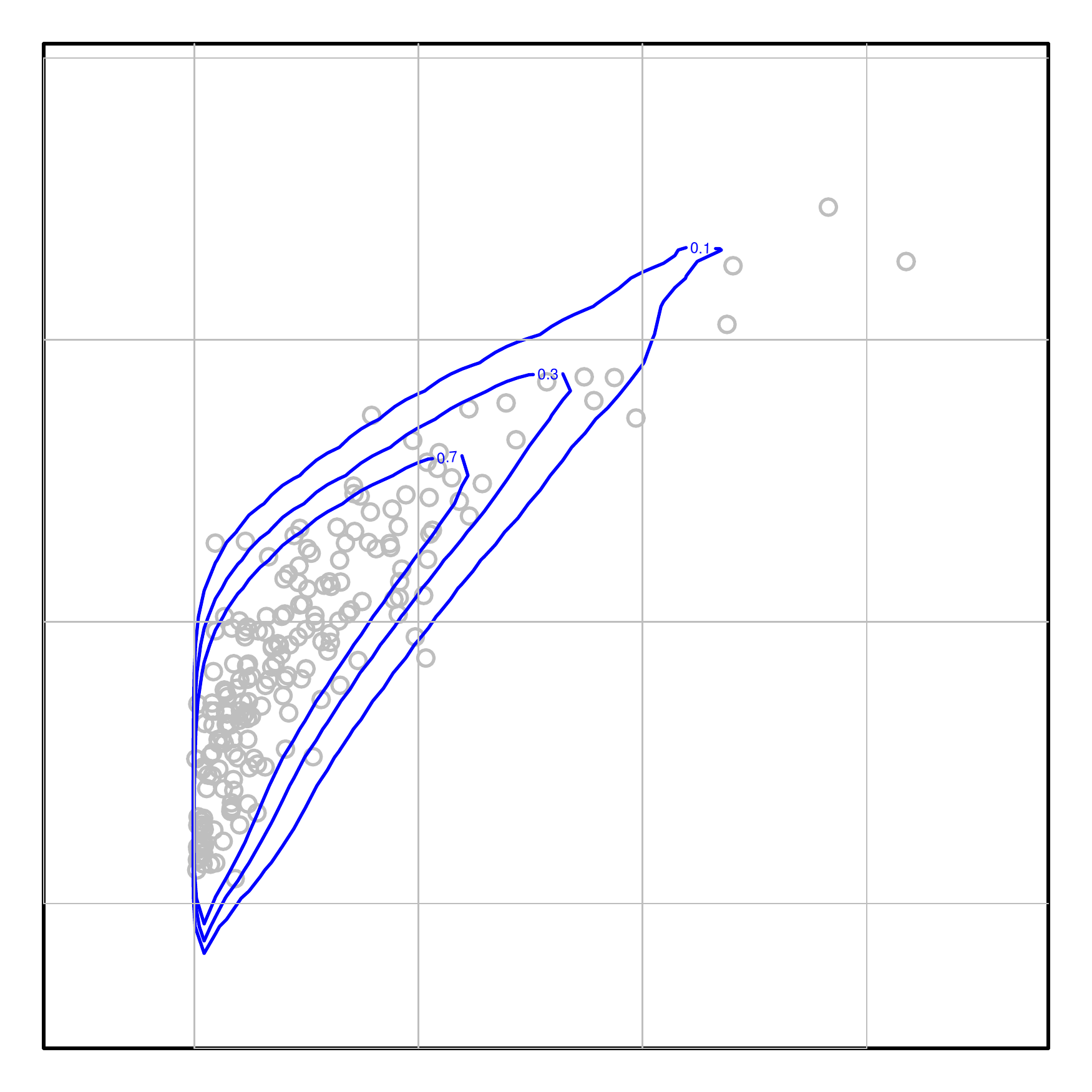} }
  \caption{Contour plots of the true
    \protect\subref{fig:subfig:n200-script4-truth}
    and the estimated joint densities with the parametric 
    \protect\subref{fig:subfig:n200-script4-para},
    nonparametric 
    \protect\subref{fig:subfig:n200-script4-nonpara},
    and semiparametric
    \protect\subref{fig:subfig:n200-script4-semi} 
    strategies. The size of the dataset is $n=200$.}
  \label{fig:n200-script4}
\end{figure}

The results for $n=200$ and marginal estimation are shown in Figure~\ref{fig:n200-script4-marginals}. Figure~\ref{fig:n200-script4-marginals}~\protect\subref{fig:subfig:n200-script4-marginal1} pictures the estimated densities with the parametric, nonparametric and semiparametric methods for the first marginal, that is, when the parametric model is well specified. We see only two lines because the parametric, semiparametric and the true densities are very similar, indicating that $\hat{\alpha} \approx 1$. Figure~\ref{fig:n200-script4-marginals}~\protect\subref{fig:subfig:n200-script4-marginal2} corresponds to the misspecified second marginal. Here this is the nonparametric and the semiparametric estimates which are nearly identical, indicating that $\hat{\alpha} \approx 0$.

Figure~\ref{fig:n200-script4} shows the estimation for the bivariate joint density. In Figure~\ref{fig:n200-script4}~\protect\subref{fig:subfig:n200-script4-para} we see that one marginal misspecification led to a poor estimation of the joint density, especially in the joint tails. Figure~\ref{fig:n200-script4}~\protect\subref{fig:subfig:n200-script4-nonpara} shows the estimated joint density with the nonparametric strategy for the marginals. Drawbacks of nonparametric estimation are easily spotted: the estimated density is multimodal and assumes positive values where it should be null. Visually, the best performance is achieved with the semiparametric strategy in Figure~\ref{fig:n200-script4}~\protect\subref{fig:subfig:n200-script4-semi}. The figures for $n=50, 100, 500$ are similar and not shown to limit the length of the paper.

\begin{figure}[t]
 \centering
 \includegraphics[width=.5\textwidth]{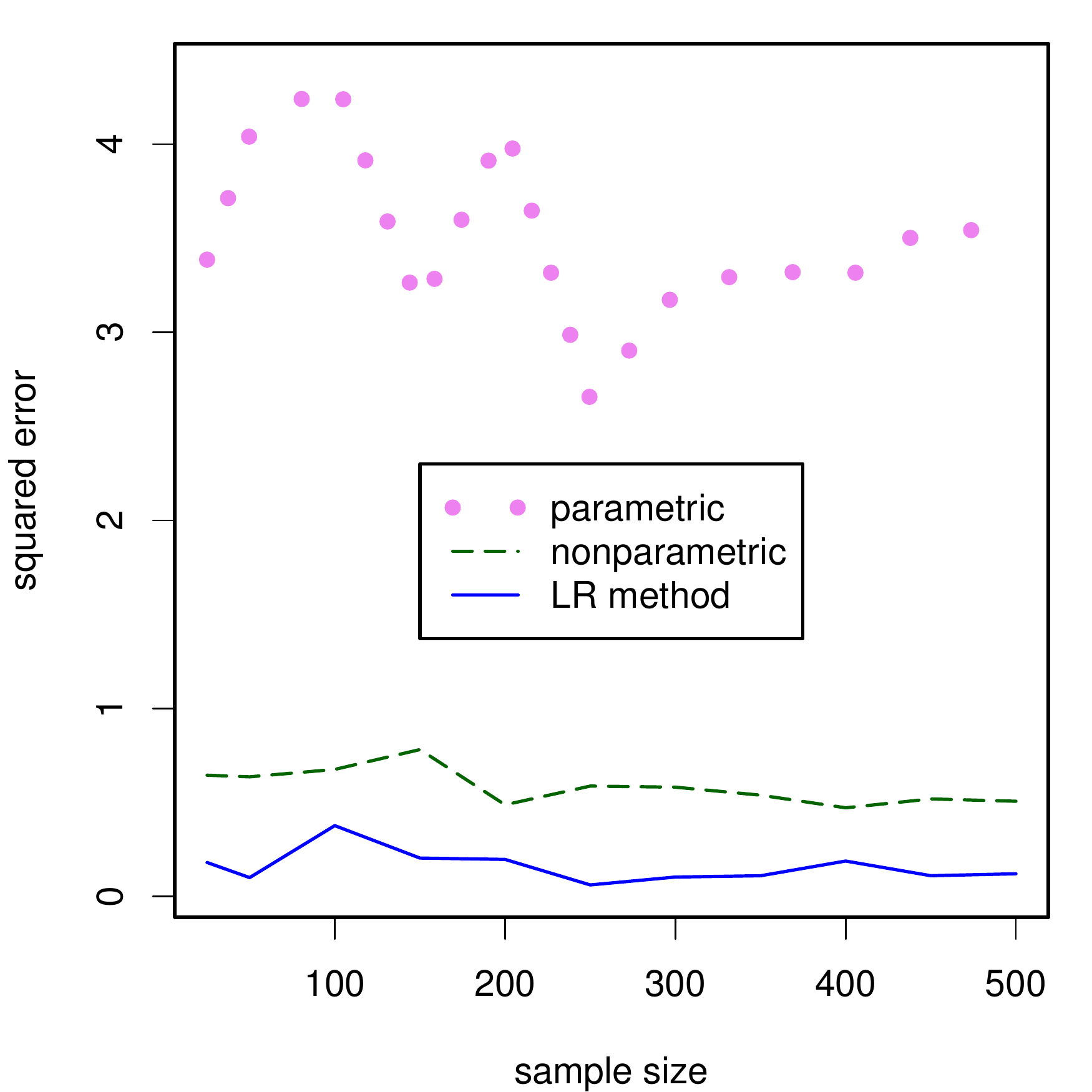}
 \caption{Squared L2-distances between the true joint density and the estimators in function of the sample size. From bottom to top, the plain blue line, green dashed line and violet dotted line are the semiparametric, the nonparametric and the parametric error curves, respectively.}
 \label{fig:errors-script4}
\end{figure}

The squared L2-distances between the true joint density and the estimators are shown in Figure~\ref{fig:errors-script4}. The semiparametric strategy performs best for all sample sizes.

\clearpage
\newpage

\begin{appendices}

\section{Proofs of the propositions}\label{app:proof_prop}

We define the mixture likelihood function $L_n : [0,1] \to \mathbb [-\infty , +\infty)$ as
\begin{align*}
L_n(\alpha)  =  \sum _ { i = 1 } ^{n }  \log \left({\alpha} f_{\hat \theta_n} (X_i) +(1-{\alpha}) \hat f _{n,i} ^{\text{LR}}  \right).
\end{align*}
The fitness coeficient $\hat \alpha_n$ in~(\ref{eq:our method}) is then defined as a maximizer of $L_n(\alpha)$ over $[0,1]$.

\subsection{Proof of Proposition \ref{prop:existence}}

The presence of $\Delta_n q(X_i)$ in $\hat f _{n,i} ^{\text{LR}}$ allows for 
$L_n(\alpha) > -\infty$ for all $\alpha \in [0,1)$. If for all $i$, $f_{\hat \theta_n} (X_i)>0$, i.e., $L_n(1)>-\infty$, then $L_n$ is continuos on $[0,1]$ and the extreme value theorem yields the existence of $\hat \alpha_n$. Else, if $L_n(1)=-\infty$, there exists $\delta>0$ such that 
$\sup_{\alpha\in[0,1-\delta]}L_n(\alpha) > \sup_{\alpha\in(1-\delta,1]}L_n(\alpha)$, meaning that the maximum is over $[0,1-\delta]$ and exists in virtue of the extreme value theorem. Whenever $f_{\hat \theta_n}(X_i)$ is not identically equal to $\hat f _{n,i} ^{\text{LR}}$ for all $i=1,\ldots,n$, the function $ L_n$ is strictly concave and so comes the unicity.

\qed

\subsection{Proof of Proposition \ref{prop:g_0tof_0}}
Let $ 0<\epsilon <1$. By assumption, there exists $\tilde A>0$, such that for all $|x|> \tilde A$, we have 
\begin{align*}
(1-\epsilon)g_0(x)  \leq f_0(x) \leq g_0(x) (1+\epsilon). 
\end{align*}
For $t >0$ small enough (i.e., taking any $t < \inf_{ |x|\leq  \tilde A}f_0(x)$ implies that $\{|x|\leq  \tilde A\} \subset \{f_0(x) >t\}$, or equivalently that $S_t^c \subset \{|x|> \tilde A\} $), it holds
\begin{align*}
\int_{S_t^c} f_0(x) \, \diff x  =  \int_{f_0(x) \leq  t } f_0(x) \, \diff x \leq  (1+\epsilon) \int_{ (1-\epsilon)  g_0(x) \leq t } g_0(x)\, \diff x.
\end{align*}
Consequently, we obtain that $\int_{S_t^c} f_0(x) \, \diff x\leq  t^{\beta} c_2(1+\epsilon) /(1-\epsilon)^{\beta}$.

Remark that 
\begin{align*}
\lambda (S_t) \leq \lambda(\{|x|\leq \tilde A\}) +  \lambda( \{|x|> \tilde A\} \cap S_t)
\leq \lambda(\{|x|\leq \tilde A\}) +  \int_{(1+\epsilon ) g_0(x)> t} dx ,
\end{align*}
which is enough to obtain the last point of (H\ref{ass:Hcvlin_f_0_St}).

Suppose that $0 < h_n\leq 1$. By enlarging $\tilde A$ (i.e., taking $\tilde A: = \tilde A +\sqrt d$), we have, for all $|x|> \tilde A$ and $u\in [-1,1]^d$,
\begin{align*}
(1-\epsilon) g_0(x+h_nu)  \leq f_0(x+h_nu)  \leq g_0(x+h_nu) (1+\epsilon).
\end{align*}
Let $b_n = \gamma (nh_n^d)^{-1/\beta}$ and $A_1 = \max(A,\tilde A)$. As soon as $|x|\leq  A_1$,
\begin{align*}
\sup_{u\in  [-1,1]^d } \frac{f_0(x+h_nu)}{f_0(x)}\leq \frac{\|f_0\|_{\mathbb R^d}}{ \inf_{|x|\leq  A_1 }{f_0(x)}}.
\end{align*}
Otherwise,
\begin{align*}
\sup_{|x| >  A_1,\, f_0(x) >  b_n } \sup_{u\in  [-1,1]^d } \frac{f_0(x+h_nu)}{f_0(x)} &\leq \frac{(1+\epsilon)}{(1-\epsilon)}  \sup_{ |x| >  A_1 ,\, (1+\epsilon)g_0(x)  >  b_n }  \sup_{u\in  [-1,1]^d } \frac{g_0(x+h_nu)}{g_0(x)}.
\end{align*}
We conclude by remarking that the previous is bounded, by $(1+\epsilon ) C_2 / (1-\epsilon) $.
\qed

\section{Proof of Theorem \ref{main:result}}\label{app:proof_main_result}

Theorem \ref{main:result} follows from the application of two high-level results, corresponding respectively to the well-specified and misspecified case. Both high-level results take place in the following general framework: given a triangular sequence of non-negative real numbers $ \xi_{n,i}, i=1,\ldots, n$, $n\geq 1$, we consider the mixture likelihood function given by
\begin{align*}
 L_n(\alpha) = \sum _ { i = 1 } ^{n }  \log \left({\alpha} f_{\hat \theta_n} (X_i) +(1-{\alpha}) \xi_{n,i}\right) .
\end{align*}
Here the sequence $(\xi_{n,i})$ is left unspecified in order to highlight the assumptions that we need on the nonparametric part. This random sequence could be the non-parametric estimator evaluated at $X_i$, i.e., $\hat f _n(X_i)$, the LOO estimate $\hat f _{n,i} ^{\text{LOO}}$ or the LR estimate $ \hat f _{n,i} ^{\text{LR}}$ with $\Delta_n >0 $. In this slightly new context, we define $\hat \alpha_n$ as
\begin{align*}
\hat \alpha_n  \in \text{argmax}_{ \alpha\in[0,1]} L_n(\alpha) .
\end{align*}
In both cases, respectively, the misspecified and well-specified case, the approach taken is similar. We compare the empirical likelihood of the mixture to the one of the parametric estimate (in the well-specified case) or the nonparametric estimate (in the misspecified case).

In the proofs below, it is convenient to introduce the normalized version of $L_n(\alpha)$, given by
\begin{align*}
\tilde {L}_n(\alpha) =  \sum _ { i = 1 } ^{n }  \log \left(\frac{ \alpha f_{\hat \theta_n,i} + (1-\alpha) \xi_{n,i}} {f_{0,i} } \right),
\end{align*}
where, for any real valued function $f$, we have introduced the short-cut notation $f_i$ for $f(X_i)$.

\subsection{Case (\ref{main:result_h0}) : the model is well-specified}

We are based on some restricted mean quadratic error
\begin{align*}
&Q_n^{\p}(S) = \sum_{i=1} ^n \left(\frac{    f_{\hat \theta_n,i} - f_{0,i }}{f_{0,i}}\right)^2 \mathrm 1_{\{ X_i \in S\}}, \qquad Q_n^{\np}(S) = \sum_{i=1} ^n  \left(\frac{  \xi_{n,i} - f_{0,i }}{f_{0,i}}\right)^2  \mathrm 1_{\{ X_i \in S\}},
\end{align*}
and some averaged linear error 
\begin{align*}
&M_n^{\p}  = \sum_{i=1} ^n \left(\frac{    f_{\hat \theta_n,i} - f_{0,i }}{f_{0,i}}\right) ,   \qquad  M_n^{\np}  = \sum_{i=1} ^n   \left(\frac{  \xi_{n,i} - f_{0,i }}{f_{0,i}}\right).
\end{align*}
The proof of the following theorem is given in Section \ref{proof:th_high_level1}.

\begin{app_theorem}\label{theorem:highlevel_condition}
 Suppose that $f_0\in \mathcal P$ and let $S\subset \mathbb R^d$ and $b>0$ be such that for all $x\in S$, $f_0(x)>b$. If the following convergences hold in probability, as $n\to \infty$,
 \begin{align}
&\| f_{\hat \theta_n} - f_{0 }\|_S  \to 0, & \max_{i=1,\ldots, n,\, :\, X_i\in S}  |\xi_{n,i} - f_{0,i}| \to 0, \label{cond:uniform_local_para+npar} \\
&\frac{Q_n^{\p}(S)}{ Q_n^{\np}(S)} \rightarrow 0 , &   \frac{|M_n^{\p}|+|M_n^{\np}| }{ Q_n^{\np}(S)} \rightarrow 0,\label{cond:high_level_rate_quad_and_lin}
\end{align} 
then, $\hat \alpha_n \rightarrow 1$ as $n\to \infty$, in probability.
\end{app_theorem}

We now verify the conditions of the previous theorem when $\xi_{n,i}$ is the LR sequence $\hat f _{n,i} ^{\text{LR}}$ and when (H\ref{ass:base_density}), (H\ref{ass:kernel}), (H\ref{ass:Hcvlin_f_0_St}), (A\ref{ass:parametric_model_consistency}), (A\ref{ass:parametric_model_assymptotic_normality}) and $nh_n^d\Delta_n \to 0$ are fulfilled. 

\paragraph{Condition (\ref{cond:uniform_local_para+npar}).}
The first convergence in (\ref{cond:uniform_local_para+npar}) holds in virtue of \eqref{eq:lip_theta_f} established in Section \ref{app:parametric}. For the second one, it holds that
\begin{align*}
\hat f _{n,i} ^{\text{LR}}    = \left(\frac{n}{n-1}\right) \left(\hat f_n(X_i) - \frac{K(0)}{nh_n^d} \right)+  \Delta_n q(X_i) .
\end{align*}
Applying the first statement of Proposition C.1 in \cite{portier+s:2015} (which is a consequence of Theorem 2.1 in \cite{gine+g:02}), we have, under (H\ref{ass:base_density}) and (H\ref{ass:kernel}), that
\begin{align*}
\|\hat f _n - f_{h_n}   \|_{\mathbb R^d} = O_{\mathbb P} \left( \sqrt{ \frac{  |\log( h_n )| }{nh_n^d }}\right)  .
\end{align*}
Together with Lemma \ref{lemma:auxiiliary_f_conv}, we obtain that $\|\hat f _n - f_{0}  \|_{\mathbb R^d}  = O_{\mathbb P} ( \sqrt{ {  |\log( h_n )| } / {nh_n^d }})  +O(h_n^2)  $.
Consequently, we get
\begin{align}
\max_{ i=1,\ldots, n}  |\hat f _{n,i} ^{\text{LR}} -f_{0,i}| &\leq  \left( \frac{n}{n-1}\right) \left(\| \hat f_{n}  -f_0 \|_{\mathbb R^d}   + \frac{K(0)}{nh_n^d}\right) +\Delta_n \|q\|_{\mathbb R^d}  + \frac{\|f_0\|_{\mathbb R^d}}{n-1}\nonumber \\
& = O_{\mathbb P} \left( \sqrt{ \frac{  |\log( h_n )| }{nh_n^d }} +h_n^2 +\Delta_n \right).\label{boud:density_unif} 
\end{align} 
The latter bound indeed goes to $0$, in probability, as $n\to \infty$.

\paragraph{Condition (\ref{cond:high_level_rate_quad_and_lin}).} We proceed as follows, with $\xi_{n,i} = \hat f _{n,i} ^{\text{LR}}$:

\begin{enumerate}[(a)]
\item By Lemma \ref{lemma:cvquad}, stated in Section \ref{app:lemma}, there exists $c>0$ such that with probability going to $1$, $h_n^d Q_n^{\np}(S)  \geq  c$. The set $S$ is chosen equal to $\{f_0(x) >b\}$ where $b>0$ is such that it is non-empty.
\item We show in Lemma \ref{lemma:cvlin} that $ h_n^d |M_n^{\np}|\to 0$ in probability.
\item In Lemma  \ref{lemma:cvquad_parametric} (resp. Lemma \ref{lemma:cvlin_param}), it is established, under (A\ref{ass:parametric_model_consistency}) and (A\ref{ass:parametric_model_assymptotic_normality}), that $Q_n^{\p}(S) = O_{\mathbb P} (1)$ (resp. $|M_n^{\p}| = O_{\mathbb P} (1)$).
\end{enumerate}
All this together implies that (\ref{cond:high_level_rate_quad_and_lin}) holds true.

\subsection{Case (\ref{main:result_h1}) : the model is misspecified}\label{sec:misspecification}

When the model is misspecified, i.e., $f_0\notin \mathcal P$, the following high-level conditions are enough to ensure the convergence in probability $\alpha_n\to 0$. These conditions are easily implied by (A\ref{ass:parametric_model_consistency}), (H\ref{ass:base_density}), (H\ref{ass:kernel}), (H\ref{ass:Hcvlin_f_0_St}) and $|\log(\Delta_n)|^{1/\beta} (\sqrt{|\log(h_n)|/nh_n^d} +h_n^2)$ as demonstrated below. The proof of this theorem is given in Section \ref{proof:th_high_level2}.

\begin{app_theorem}\label{theorem:high_level_2}
Suppose that $f_0\notin \mathcal P$ and $\mathbb E [|\log(f_{0,1})|] < \infty$, that the class $\mathcal P $ is Glivenko-Cantelli (i.e., (\ref{eq:GC}) holds) with $\Theta$ compact, that the envelop $F_{\Theta} $ is such that $\mathbb E [\log(F_{\Theta,1})]<+\infty$, and that for every $x\in \mathbb R^d $, $\theta\mapsto f_\theta(x)$ is a continuous function defined on $\Theta$. Suppose that there exists $\beta\in(0,1] $ and $c>0$ such that, as $t\to 0$,  $\int_{S_t^c} f_0(x) \leq ct^{\beta}$, and $q:\mathbb R^d\to \mathbb R^+$ such that $\mathbb E[|\log(q(X_1))|]<\infty$, $ \xi_{n,1} \geq \Delta_n q_1 $ a.s., and
\begin{align*}
&  |\log(\Delta_n)|^{1/\beta}  \max_{i=1,\ldots,n}   | \xi_{n,i}   -  f_0(X_i) | \rightarrow 0,
\end{align*}
then, $\hat \alpha_n \to 0$, as $n\to \infty$, in probability.
\end{app_theorem}

We already argued that (\ref{eq:GC}) is implied by (A\ref{ass:parametric_model_consistency}). The continuity of $f_\theta$ is deduced from the continuity of $\log(f_\theta)$ provided by (A\ref{ass:parametric_model_consistency}).  The bound given in (\ref{boud:density_unif}) together with $|\log(\Delta_n)|^{1/\beta} (\sqrt{|\log(h_n)|/nh_n^d} +h_n^2)\to 0$ implies the stated convergence with $\xi_{n,i} = \hat f _{n,i} ^{\text{LR}}$.

\subsection{Proofs of the high-level Theorems}
\subsubsection{Proof of Theorem \ref{theorem:highlevel_condition}}\label{proof:th_high_level1}

Because $f_0 \in \mathcal P$, it holds that $f_{\hat \theta_n,i}>0$ for all $i=1,\ldots, n$, which guarantees the existence of a maximizer $\hat \alpha_n$ (as explained in the proof of Proposition \ref{prop:existence}). 
By definition of $\hat \theta_n$,  $ \sum _ { i = 1 } ^{n }  \log \left( f_{0,i} \right)\leq  \sum _ { i = 1 } ^{n }  \log \left( f_{\hat \theta_n,i} \right)$. Consequently, $\max_{\alpha\in [0,1]} \tilde {L}_n(\alpha)\geq 0 $ and for every $\epsilon>0$, the event $\max_{\alpha\in [0,1-\epsilon]} \tilde {L}_n(\alpha)<0$ implies that $\hat \alpha_n>1- \epsilon$. Thus, let $\epsilon\in (0,1)$, the proof will be completed by showing that with probability going to $1$,
\begin{align*}
\sup_{\alpha\in [0, 1-\epsilon]}  \tilde {L}_n(\alpha) <0.
\end{align*}
A useful notation in the following is
\begin{align*}
\hat x_{i,n} =1 +  \frac{  \alpha (  f_{\hat \theta_n,i} - f_{0,i })}{f_{0,i}} +\frac{  (1-\alpha) ( \xi_{n,i} - f_{0,i })}{f_{0,i}} .
\end{align*} 
A useful technical detail is there exists a sequence $\epsilon_n\to 0$ such that the event 
\begin{align*}
\{  \max_{i=1\ldots,n\, : \,X_i\in S } {f_{0,i}} |\hat x_{i,n}-1|\leq \epsilon_n\}
\end{align*} 
has probability going to $1$ as $n\to \infty$. This is a consequence of (\ref{cond:uniform_local_para+npar}). As we are establishing a result in probability, we can further suppose that this event is realized.

A key step in our approach is the following inequality, reminiscent of the Taylor development of the logarithm around $1$,
\begin{align*}
\log(x) - (x-1) \leq \left\{ \begin{array}{lr}
 -\frac{1}{4}  (x-1)^2  &\text{if } 1/2 < x < 3/2 \\
 0 &\text{else }   
\end{array}\right. ,
\end{align*}
which might be derived by studying the concerned function. This kind of inequality is commonly used for studying likelyhood methods \citep{grenander:1981,murphy:1994}. Applied to $\hat x_{i,n}$, it gives
\begin{align*}
\tilde {L}_n (\alpha ) -(\alpha M_n^{\p} + (1-\alpha) M_n^{\np}  )&= \sum_{i=1} ^n (  \log\left(\hat x_{i,n} \right)   -(\hat x_{i,n}  -1) )\\
&\leq -\frac{1}{4}\, \sum_{i=1} ^n (\hat x_{i,n}-1)^2\mathrm 1_{\{|\hat x_{i,n}-1| < 1/2 \}}  \\
&\leq -\frac{1}{4}\, \sum_{i=1} ^n (\hat x_{i,n}-1)^2\mathrm 1_{\{X_i\in S, \, |\hat x_{i,n}-1| < 1/2\}} .
\end{align*}
Note that whenever $X_i\in S$, because it holds ${f_{0,i}} |\hat x_{i,n}-1|\leq \epsilon_n$, we have (for $n$ small enough) that $|\hat x_{i,n}-1|<1/2$. This means that, for all $i=1,\ldots, n$, $   \mathrm 1_{\{X_i\in S\} } \leq \mathrm 1_{\{|\hat x_{i,n}-1| < 1/2 \}}$, and it follows
\begin{align*}
\tilde {L}_n (\alpha ) -(\alpha M_n^{\p} + (1-\alpha) M_n^{\np}  )  &\leq -\frac{1}{4}\, \sum_{i=1} ^n (\hat x_{i,n}-1)^2\mathrm 1_{\{X_i \in S\}}  \\
& = - \frac{1}{4}\left\{  (1-\alpha)^2 Q_n^{\np}(S)  + \alpha^2 {Q_n^{\p}(S)} +2  \alpha(1-\alpha)   {U_n}\right\} \\
&\leq  - \frac{1}{4} (1-\alpha)^2 Q_n^{\np}(S) \left\{  1 -  \frac{ 2\alpha  {|U_n|}  }{ (1-\alpha) Q_n^{\np}(S)}  \right\},
\end{align*}
where
\begin{align*}
&U_n = \sum_{i=1} ^n \frac{   (  f_{\hat \theta_n,i} - f_{0,i }) ( \xi_{n,i} - f_{0,i })}{f_{0,i}^2}  \mathrm 1_{\{ X_i \in S\}} .
\end{align*}
Bounding the right-hand side with respect to $\alpha\in [0,1-\epsilon]$ gives
\begin{align*}
\sup_{\alpha\in [0, 1-\epsilon]} \{ \tilde {L}_n (\alpha ) -(\alpha M_n^{\p} + (1-\alpha) M_n^{\np}  ) \} &\leq  - \frac{1}{4} \epsilon ^2  Q_n^{\np}(S) \left(   1   -  2 \epsilon ^{-1} \frac{|U_n|}{Q_n^{\np}(S)}\right).
\end{align*}
By assumption, we have that $ {Q_n^{\p}(S) }/{Q_n^{\np}(S)}\to 0$ in probability. From the Cauchy-Schwartz inequality we get that $|U_n|\leq \sqrt{Q_n^{\np}(S)Q_n^{\p}(S)}$, leading to $ {|U_n| }/{Q_n^{\np}(S)}\to 0$, in probability. Consequently, we obtain that
\begin{align*}
 \sup_{\alpha\in [0, 1-\epsilon]}  \tilde {L}_n(\alpha)  & \leq -\frac{1}{4} \epsilon ^2  Q_n^{\np}(S) \left(   1   -  2 \epsilon ^{-1} \frac{|U_n|}{Q_n^{\np}(S)}  - 4 \frac{| M_n^{\p} |+ | M_n^{\np} | }{ \epsilon ^2 Q_n^{\np}(S)  }    \right) .
\end{align*}
The term between brackets goes to $1$, in probability, implying that for every $\delta>0$, with probability going to $1$,
\begin{align*}
   \sup_{\alpha\in [0, 1-\epsilon]} \tilde {L}_n(\alpha, \hat \theta)  & \leq   -\frac{1}{4} \epsilon ^2  Q_n^{\np}(S) ( 1 - \delta) .
\end{align*}
Hence it remains to note that, by (\ref{cond:high_level_rate_quad_and_lin}), with probability going to $1$, $ Q_n^{\np}(S)>0$.

\qed

\subsubsection{Proof of Theorem \ref{theorem:high_level_2}}\label{proof:th_high_level2}

Note that $\hat \alpha _ n \in \argmax_{\alpha\in [0,1]} \tilde L_n(\alpha) $ exists because $\xi_{n,i}>0$ for all $i$, as explained in the proof of Proposition \ref{prop:existence}. Let $\epsilon>0$. The proof requires to show that with probability going to $1$, $\hat \alpha_n <\epsilon$. This event is realized as soon as $\max_{\alpha\in [\epsilon, 1]} \tilde {L}_{n} (\alpha) < \tilde {L}_n (0)$. We analyse both terms separately. First we show that 
 \begin{align*}
 \tilde {L}_n (0)\to 0,
 \end{align*}
in probability, and then that there exists $\delta>0$ such that, with probability going to $1$,
 \begin{align}\label{eq:proof_fact2}
\sup_{\alpha\in [\epsilon, 1]}\tilde {L}_n (\alpha) \leq -\delta.
 \end{align}
Let $\eta>0$, $b_n = (\eta / |\log(\Delta_n)|)^{1/\beta} $ and $c_n = \max_{i=1,\ldots,n} |\xi_{n,i} - f_{0,i} |$. We assume further that $b_n+c_n<1$ and $\Delta_n<1$.
 We have 
\begin{align*}
&| \tilde L_{n} (0)| \\
&\leq \left|  n ^{-1}  \sum_{i=1} ^n  \log\left(\frac{  \xi_{n,i}}{f_{0,i}} \right) \mathrm 1_{\{  f_{0,i} >  b_n \}} \right|  +  n ^{-1}  \sum_{i=1} ^n  \left| \log\left(\frac{  \xi_{n,i}}{f_{0,i}} \right)\right| \mathrm 1_{\{ f_{0,i} \leq b_n\}}\\
& \leq  \left|  n ^{-1}  \sum_{i=1} ^n  \log\left(\frac{  \xi_{n,i}}{f_{0,i}} \right) \mathrm 1_{\{  f_{0,i} >  b_n \}} \right|  +n ^{-1}  \sum_{i=1} ^n  \left( \left| \log(\Delta_nq_i )\right|    \mathrm 1_{\{ f_{0,i} \leq  b_n\}}+  \left| \log(f_{0,i})\right| \mathrm 1_{\{ f_{0,i} \leq b_n\}}\right)\\
&\leq  \left|  n ^{-1}  \sum_{i=1} ^n  \log\left(\frac{  \xi_{n,i}}{f_{0,i}} \right) \mathrm 1_{\{  f_{0,i} >  b_n \}} \right|  + \left| \log(\Delta_n)\right|  n ^{-1}  \sum_{i=1} ^n     \mathrm 1_{\{ f_{0,i} \leq  b_n\}}\\
&\hspace{3cm}+ n ^{-1}  \sum_{i=1} ^n  (\left| \log(q_{i})\right|  + \left| \log(f_{0,i})\right| )\mathrm 1_{\{ f_{0,i} \leq b_n\}}
\end{align*}
The expectation of the term in the middle is bounded by $ | \log(\Delta_n )|  \mathbb P ( f_{0} (X_1 ) \leq b_n  )  $
of order $ | \log(\Delta_n )| b_n^{\beta} = \eta $, by assumption. The corresponding term goes to $0$, as $\eta $ is arbitrarily small. The expectation of the term in the right is smaller than $\mathbb E[ (|\log(q_1)|+  | \log(f_{0,1})|) \mathrm 1_{\{ f_{0,i} \leq b_n\}}]$, which goes to $0$ because $|\log(q_{1})|$ and $|\log(f_{0,1})| $ are integrable. Hence it remains to obtain that the term in the left goes to $0$. The mean-value theorem gives
\begin{align*}
\left| n ^{-1}  \sum_{i=1} ^n  \log\left(\frac{  \xi_{n,i}}{f_{0,i}} \right) \mathrm 1_{\{  f_{0,i} >b_n \}} \right| &\leq   \frac{\max_{i=1,\ldots,n} |\xi_{n,i} - f_{0,i} |  }{\inf_{i=1,\ldots, n\,:\, f_{0,i}>b_n} \inf_{t\in [0,1]}  t\xi_{n,i}+ (1-t) {f_{0,i} }  } \\
&\leq \frac{c_n}{b_n - c_n  }\to 0.
\end{align*}
Now we establish (\ref{eq:proof_fact2}) by obtaining one-sided inequalities. Take $b_n = (1 / |\log(\Delta_n)|)^{1/\beta} $, suppose that $b_n+c_n<1$, and use the monotonicity of the logarithm, to get that
\begin{align*}
 n ^{-1} \sum_{i=1} ^n \log\left(\frac{ \alpha f_{\hat \theta_n,i} + (1-\alpha) \xi_{n,i} }{f_{0,i}} \right)    \mathrm 1_{\{  f_{0,i} \leq b_n   \}} &\leq    n ^{-1} \sum_{i=1} ^n  \log\left(\frac{  F_{\Theta , i} +  b_n  + c_n  }{ f_{0,i} } \right)  \mathrm 1_{\{  f_{0,i} \leq  b_n   \}}\\
&\leq    n ^{-1} \sum_{i=1} ^n | \log\left(\frac{  F_{\Theta , i} + 1 }{ f_{0,i} } \right)|  \mathrm 1_{\{  f_{0,i} \leq  b_n   \}}.
\end{align*}
Taking the expectation, we find a bound in $ \mathbb E[|\log\left( \frac{  F_{\Theta , 1} +  1  }{ f_{0,1} } \right) | \mathrm 1_{\{  f_{0,1} \leq  b_n   \}}] $ which goes to $0$ as $n\to \infty$ in virtue of the Lebesgue dominated convergence theorem. Then, taking $0<\eta<1$, it holds that
\begin{align*}
\tilde L_{n} (\alpha) &\leq n ^{-1} \sum_{i=1} ^n \log\left(\frac{ \alpha f_{\hat \theta_n,i} + (1-\alpha) \xi_{n,i} }{f_{0,i}} \right)    \mathrm 1_{\{  f_{0,i} > b_n \}} + o_p(1) \\
&\leq  n ^{-1} \sum_{i=1} ^n \log\left(\frac{ \alpha (f_{\hat \theta_n,i}+\eta f_{0,i}) + (1-\alpha) \xi_{n,i} }{f_{0,i}} \right)    \mathrm 1_{\{  f_{0,i} > b_n \}} + o_p(1)
\end{align*}
The first term in the right-hand side is decomposed according to
\begin{align*}
n ^{-1} \sum_{i=1} ^n &\log\left(\frac{ \alpha (f_{\hat \theta_n,i}+\eta f_{0,i}) + (1-\alpha) \xi_{n,i} }{ \alpha (f_{\hat \theta_n,i}+\eta f_{0,i}) + (1-\alpha) f_{0,i}} \right)    \mathrm 1_{\{  f_{0,i} > b_n \}} \\
& +n ^{-1} \sum_{i=1} ^n \log\left(\frac{ \alpha (f_{\hat \theta_n,i}+\eta f_{0,i}) + (1-\alpha) f_{0,i} }{ f_{0,i}} \right)    \mathrm 1_{\{  f_{0,i} > b_n \}} .
\end{align*}
By the mean value theorem, the term on the left is bounded by
\begin{align*}
\frac{  (1-\alpha) \max_{i=1,\ldots, n} |\xi_{n,i}  - f_{0,i}|}{ (\eta b_n) \wedge (b_n - c_n)},
\end{align*}
which goes to $0$, by assumption. For the term on the right, notice that $\{   \alpha (f_{\theta}+\eta f_{0}) + (1-\alpha) f_{0} \, :\, \alpha\in [\epsilon, 1],\, \theta \in \Theta \} $ is Glivenko-Cantelli with envelop $ F_{\Theta} + 2 f_{0}$. Then applying Theorem 3 in \cite{vandervaart+w:2000}, the class formed by $\log( \alpha (f_{\theta}+\eta f_{0})  + (1-\alpha) f_{0} ) $ is still Glivenko-Cantelli. Since for all $\theta\in \Theta$, $\alpha\in [\epsilon,1]$,
\begin{align*}
\log( \epsilon \eta f_{0}  ) \leq \log( \alpha (f_{\theta}+\eta f_{0})  + (1-\alpha) f_{0} )\leq \log(  F_{\Theta}+2 f_{0} +1),
\end{align*}
 the function $|\log( \epsilon \eta f_{0}  )| + \log(  F_{\Theta}+2f_0 +1 )$ is an integrable envelop. Using again Theorem 3 in \cite{vandervaart+w:2000}, the class formed by $\log( \alpha (f_{\theta}+\eta f_{0})  + (1-\alpha) f_{0} )\mathrm 1_{\{ f_{0} >b\} } $, $\theta\in \Theta$, $\alpha\in [\epsilon,1]$, $0<b<1$, is still Glivenko-Cantelli with the same envelop. This implies that
\begin{align*}
\sup_{ \alpha \in [\epsilon , 1],\,  \theta\in \Theta} \left| n ^{-1} \sum_{i=1} ^n \log\left(\frac{  \alpha (f_{\theta,i}+\eta f_{0,i})+ (1-\alpha) f_{0,i} }{ f_{0,i}} \right) \mathrm {1} _ { \{f_{0,i}>b_n\} } \right.\\
\left.-\mathbb E\left[\log\left(\frac{  \alpha (f_{\theta,1}+\eta f_{0,1}) + (1-\alpha) f_{0,1} }{ f_{0,1}} \right)\mathrm {1} _ { \{f_{0,1}>b_n\} }   \right] \right| \rightarrow 0.
\end{align*}
The integrability of the envelop and the fact that $b_n\to 0$ implies that
\begin{align*}
\sup_{ \alpha \in [\epsilon , 1],\,  \theta\in \Theta} \mathbb E\left[\log\left(\frac{  \alpha (f_{\theta,1}+\eta f_{0,1}) + (1-\alpha) f_{0,1} }{ f_{0,1}} \right)\mathrm {1} _ { \{f_{0,1}\leq b_n\} }   \right] \to 0.
\end{align*} 
It remains to use the inequality $\log(x)\leq 2(\sqrt x - 1)$ to obtain that
\begin{align*}
\sup_{ \alpha \in [\epsilon , 1],\,  \theta\in \Theta} \mathbb E\left[\log\left(\frac{ m_{\theta,\alpha,\eta} }{ f_{0,1}} \right)   \right]&\leq \sup_{ \alpha \in [\epsilon , 1],\,  \theta\in \Theta} 2 \int (\sqrt{ m_{\theta,\alpha,\eta}f_0 } - f_0) \,\diff \lambda
\end{align*}
where $m_{\theta,\alpha,\eta}=  \alpha (f_{\theta}+\eta f_{0}) + (1-\alpha) f_{0}$. Since
\begin{align*}
\sup_{ \alpha \in [\epsilon , 1],\,  \theta\in \Theta} | \int \sqrt{ (m_{\theta,\alpha,\eta} - m_{\theta,\alpha,0})  f_0 }  \,\diff \lambda | = \sqrt \eta ,
\end{align*}
we get, using that $\sqrt {a+b}\leq \sqrt{a} + \sqrt b$, $a\geq 0$, $b\geq 0$,
\begin{align*}
\sup_{ \alpha \in [\epsilon , 1],\,  \theta\in \Theta} \mathbb E\left[\log\left(\frac{ m_{\theta,\alpha,\eta} }{ f_{0,1}} \right)   \right]&\leq 2\sqrt \eta + 2\sup_{ \alpha \in [\epsilon , 1],\,  \theta\in \Theta}  \int (\sqrt{  m_{\theta,\alpha,0}  f_0 }  - f_0) \,\diff \lambda \\
&= 2\sqrt \eta - \inf_{ \alpha \in [\epsilon , 1],\,  \theta\in \Theta}  \int (\sqrt{  m_{\theta,\alpha,0} }-\sqrt{  f_0 } )^2 \,\diff \lambda.
\end{align*}
Using standard results about the Hellinger distance~\cite{pollard:2000} (chapter 3) we obtain
\begin{align*}
\sup_{ \alpha \in [\epsilon , 1],\,  \theta\in \Theta} \mathbb E\left[\log\left(\frac{ m_{\theta,\alpha,\eta} }{ f_{0,1}} \right)   \right]
&\leq 2\sqrt \eta - (1/4) \inf_{ \alpha \in [\epsilon , 1],\,  \theta\in \Theta}  \alpha^2\left( \int |f_\theta - f_0| \, \diff \lambda\right)^2\\
&\leq 2\sqrt \eta - (\epsilon ^2/ 4) \left(  \inf_{  \theta\in \Theta}  \int |f_\theta - f_0| \, \diff \lambda\right)^2.
\end{align*}
Since $f_0\notin \mathcal P$ and by the continuity assumption on $f_\theta$, it holds that $\inf_{\theta\in \Theta } \int |f_\theta-f_0| \diff \lambda >0$.
Then, as $\eta$ is arbitrary, the proof of (\ref{eq:proof_fact2}) is complete.

\qed

\subsection{Linear and quadratic error of parametric and nonparametric estimate}\label{app:lemma}

Important tools for dealing with the terms involving $\hat f _{n,i} ^{\text{LR}}$ are coming from $U$-statistic theory. We call $U$-statistic of order $p$ with kernel $w:\mathbb R^p \rightarrow \mathbb R$, any quantity of the kind
\begin{align*}
\sum_{i_1,\ldots, i_p \in D} w(X_{i_1},\ldots, X_{i_p}),
\end{align*}
where the summation is taken over the subset $D$ formed by the $(i_1,\ldots, i_p) \in \{1,\ldots, n\}^p$ such that $i_k\neq i_\ell $, $\forall k\neq \ell$. The number of terms in the summation is then  $n(n-1)\ldots (n-p+1)$. When the kernel $w$ is such that, for every $k\in\{1,\ldots, p\}$, $\mathbb E[w(X_{1},\ldots, X_{p})\mid X_{1},\ldots X_{k-1},X_{k+1},\ldots, X_{p} ]=0$, it is called a degenerate $U$-statistic.  In the proofs, we shall rely on the so-called Hajek decomposition \citep[Lemma 11.11]{vandervaart:1998}. 

To establish the two following lemmas, Lemma \ref{lemma:cvquad} and Lemma \ref{lemma:cvlin}, we are based on (H\ref{ass:base_density}), (H\ref{ass:kernel}) and (H\ref{ass:Hcvlin_f_0_St}). One might note that the expressions (a) or (b) in (H\ref{ass:kernel}) on the kernel are not used in any of these lemmas.

\begin{app_lemma}\label{lemma:cvquad}
Under assumptions (H\ref{ass:base_density}), (H\ref{ass:kernel}) and (H\ref{ass:Hcvlin_f_0_St}), if $n h_n^{d}\Delta_n \to 0$, for any $\delta>0$ and any set $S\subset \mathbb R^d $ such that $\inf_{x\in S}f_0(x)>b$, we have with probability going to $1$,
\begin{align*}
& {h_n^d} \, \sum _ { i = 1 } ^{n }  \left(\frac{ \hat f _{n,i} ^{\text{LR}} - f_{0,i}}{f_{0,i}} \right)^2  \mathrm 1_{\{X_i \in  S \}}  \geq (1-\delta)  v_K\lambda (S) .
\end{align*} 
where $v_K = \int K(u) ^2 \, \diff u$.
\end{app_lemma}

\begin{proof}
 Note that
\begin{align*}
 \mathbb E  \left[ \hat f _{n,i} ^{\text{LR}}  \mid  X_i \right]  = (n-1)^{-1} \sum_{j \neq  i} ^{n} \mathbb E\left[ { h_n^{-d} K\left(\frac{X_i - X_j}{h_n} \right) }   \mid  X_i \right] +\Delta_nq_i = f_{h_n,i}+\Delta_nq_i .
\end{align*}
The proof follows from the decomposition
\begin{align*}
\sum _ { i = 1 } ^{n }\left(\frac{ \hat f _{n,i} ^{\text{LR}}- f_{0,i}}{f_{0,i}} \right) ^2\mathrm 1_{\{X_i \in  S \}}   = A_n+B_n+2C_n,
\end{align*}
where 
\begin{align*}
&A_n =  \sum _ { i = 1 } ^{n } \left(\frac{ \hat f _{n,i} ^{\text{LR}} -  \mathbb E  \left[ \hat f _{n,i} ^{\text{LR}}   \mid  X_i \right]}{f_{0,i}} \right) ^2\mathrm 1_{\{X_i \in  S \}} ,  \\
&B_n = \sum _ { i = 1 } ^{n } \left(\frac{ f_{h_n,i}+\Delta_n q_i  - f_{0,i}}{f_{0,i}} \right) ^2  \mathrm 1_{\{X_i \in  S \}},\\
&C_n =  \sum _ { i = 1 } ^{n } \frac{ \left(\hat f _{n,i} ^{\text{LR}} -  \mathbb E  \left[ \hat f _{n,i} ^{\text{LR}}  \mid  X_i \right] \right)(f_{h_n,i}+\Delta_nq_i - f_{0,i})}{f_{0,i} ^2  } \mathrm 1_{\{X_i \in  S \}} .
\end{align*}
We will show that $ h_n^dA_n \to v_K\lambda (S)$, in probability and that $h_n^d C_n \to 0$, in probability. This will be enough as $B_n\geq 0$, almost surely. 

\paragraph{Proof that $ h_n^dA_n \to v_K\lambda (S)$ in probability.}
Introduce the notation, for any $h>0$,
\begin{align*}
&a_h(x,y) = \frac{ K_h(x-y)  - f_{h}(x)}{f_{0}(x)  } ,\\
&u_{h}(x,y,z) = a_h(x,y) a_h(x,z) \mathrm 1_{\{x \in  S \}}.
\end{align*}
Developing, we find
\begin{align*}
&A_n =  (n-1)^{-2} \sum _ { i = 1 } ^{n }  \sum _ { j\neq  i} ^{n } \sum _ { k \neq  i } ^{n }u_{h_n}(i,j,k) , 
\end{align*}
where $u_{h_n}(i,j,k) $ is as short-cut for $u_{h_n}(X_i,X_j,X_k) $. We treat $A_n$ relying on the Hajek projection of $U$-statistics. Up to a centering term, $\mathbb E \left[ u_{h_n}(i,j,k) \mid X_j, X_k\right] $, the $U$-statistic $A_n$ is a degenerate $U$-statistic. In the following we voluntary introduce this centering term in the summation to handle separately a degenerate U-statistic and another summation with less indices. By introducing, for any $h>0$,
\begin{align*}
&v_{h}(j,k) = \mathbb E \left[ u_{h}(i,j,k) \mid X_j, X_k\right]  , \\
& w_{h}(i,j,k) = u_{h}(i,j,k)- v_{h}(j,k) ,
\end{align*}
we obtain
\begin{align}
\nonumber A_n &= (n-1)^{-2} \sum _ { i = 1 } ^{n }  \sum _ { j\neq  i} ^{n } \sum _ { k \neq  i } ^{n }w_{h_n}(i,j,k)  + (n-1)^{-2}  \sum _ { i = 1 } ^{n } \sum _ { j\neq  i} ^{n } \sum _ { k \neq  i } ^{n } v_{h_n}(j,k) \\
\nonumber &= (n-1)^{-2} \sum _ { i = 1 } ^{n }  \sum _ { j\neq  i} ^{n } \sum _ { k \neq  i, \, k\neq j} ^{n }w_{h_n}(i,j,k)  
+(n-1)^{-2} \sum _ { i = 1 } ^{n }  \sum _ { j\neq  i} ^{n }  w_{h_n}(i,j,j)  \\
\label{eq:A_n}  &\qquad +(n-1)^{-2}  \sum _ { i = 1 } ^{n } \sum _ { j\neq  i} ^{n } \sum _ { k \neq  i } ^{n } v_{h_n}(j,k) .
\end{align}

\noindent\textit{Treatment of the first term in (\ref{eq:A_n}).}
Note that $w_{h_n}(i,j,k)$ defines a degenerate $U$-statistic, i.e., 
\begin{align*}
\mathbb  E[ w_{h_n}(i,j,k)\mid X_i,X_j] =\mathbb  E[ w_{h_n}(i,j,k)\mid X_i,X_k] =\mathbb  E[ w_{h_n}(i,j,k)\mid X_j,X_k] =0.
\end{align*}
Note that
\begin{align*}
\sum _ { i = 1 } ^{n }  \sum _ { j\neq  i} ^{n } \sum _ { k \neq  i, \, k\neq j} ^{n }w_{h_n}(i,j,k)   = \sum _ { i = 1 } ^{n }  \sum _ { j>i} ^{n } \sum _ {  k>j } ^{n }\overline{w}_{h_n}(i,j,k) , 
\end{align*}
 where $\overline{w}_{h}$ is the symmetrized version of $w_h$, i.e., for any triplet $(x_1,x_2,x_3)$ of $\overline{w}_h(x_1,x_2,x_3) = \sum_{\sigma} w_h (x_{\sigma(1)},x_{\sigma(2)},x_{\sigma(3)})$ where the sum is over all the $3!$ possible permutations of the set $\{1,2,3\}$. Using that the $U$-statistic with kernel  $\overline{w}_{h_n}$ is degenerate, some algebra gives that 
\begin{align*}
\mathbb E\left[\left( (n-1)^{-2} \sum _ { i = 1 } ^{n }  \sum _ { j\neq  i} ^{n } \sum _ { k \neq  i, \, k\neq j} ^{n }w_{h_n}(i,j,k)  \right)^2  \right]  &=  (n-1)^{-4} \sum _ { i = 1 } ^{n }  \sum _ { j>i} ^{n } \sum _ {  k>j } ^{n } \mathbb E[\overline{w}_{h_n}(i,j,k)^2 ] \\
&= O(n^{-1}) \mathbb E[\overline{w}_{h_n}(1,2,3)^2 ] .
\end{align*} 
 We have, using Minkowski's inequality and the definition of the conditional expectation, that
\begin{align*}
\sqrt {E[\overline{w}_{h_n}(1,2,3)^2 ]}&\leq 3!   \sqrt {E[{w}_{h_n}(1,2,3)^2 ]}\leq 3!   \sqrt {E[{u}_{h_n}(1,2,3)^2 ]}.
\end{align*}
Consequently, in virtue of (\ref{ineq:order2moments1}) in Lemma \ref{lemma:auxiliary_bound_kernel}, we have shown that
\begin{align*}
\mathbb E\left[\left(  (n-1)^{-2} \sum _ { i = 1 } ^{n }  \sum _ { j\neq  i} ^{n } \sum _ { k \neq  i, \, k\neq j} ^{n }w_{h_n}(i,j,k)  \right)^2  \right] & =  O(n^{-1} h_n^{-2d} ) .
\end{align*} 
The previous rate, multiplied by $h_n^{2d}$, goes to $0$, hence, this term is negligible. 

\noindent\textit{Treatment of the second term in (\ref{eq:A_n}).}
We continue the study of $A_n$ by considering
\begin{align*}
&(n-1)^{-2} \sum _ { i = 1 } ^{n }  \sum _ { j\neq  i} ^{n }  w_{h_n}(i,j,j) \\
&\qquad = (n-1)^{-2} \sum _ { i = 1 } ^{n }  \sum _ { j\neq  i} ^{n }  (w_{h_n}(i,j,j) - \mathbb E [ u_{h_n}(i,j,j)\mid X_i]+ \mathbb E [ u_{h_n}(1,2,2)] ) \\
&\qquad\qquad +(n-1)^{-1} \sum _ { i = 1 } ^{n }  (\mathbb E [ u_{h_n}(i,j,j)\mid X_i] - \mathbb E [ u_{h_n}(1,2,2)] ) .
\end{align*}
The first term is a degenerate $U$-statistic of order $2$ whose order $2$ moments satisfy
\begin{align*}
&\mathbb E\left[\left(  (n-1)^{-2} \sum _ { i = 1 } ^{n }  \sum _ { j\neq  i} ^{n }  (w_{h_n}(i,j,j) - \mathbb E [ u_{h_n}(i,j,j)\mid X_i]+ \mathbb E [ u_{h_n}(1,2,2)] ) \right)^2  \right]\\
&= O(n^{-2}) \mathbb E [u_{h_n}(1,2,2)^2 ] = O(n^{-2}h_n^{-3d}).
\end{align*}
This is obtained by following exactly the same lines as in the treatment of the $U$-statistic $w_n$ and using (\ref{ineq:order2moments_bis_1}) in Lemma \ref{lemma:auxiliary_bound_kernel}. As $n^{-2}h_n^{-3d} \times h_n^{2d} \to  0$, the previous term is negligible. The second term is a sum of centred independent random variables with variance smaller than, in virtue of (\ref{ineq:order2moments0}) in Lemma \ref{lemma:auxiliary_bound_kernel},
\begin{align*}
n(n-1)^{-2}   \mathbb E \left[ \mathbb E [ u _{h_n}(1,2,2)\mid X_1]^2 \right] &= O(n^{-1} h_n^{-2d} ) .
\end{align*}
This is the same rate as the rate obtained for the (negligible) $U$-statistic of order $3$ with kernel $w_n$. 

\noindent\textit{Treatment of the third term in (\ref{eq:A_n}).}
The study of $A_n$ continues by considering
\begin{align*}
&(n-1)^{-2}  \sum _ { i = 1 } ^{n } \sum _ { j\neq  i} ^{n } \sum _ { k \neq  i } ^{n } v_{h_n}(j,k) \\
&= (n-1)^{-2}  \sum _ { i = 1 } ^{n } \sum _ { j\neq  i} ^{n }  v_{h_n}(j,j)
 + (n-1)^{-2}  \sum _ { i = 1 } ^{n } \sum _ { j\neq  i} ^{n } \sum _ { k \neq  i,\, k\neq  j } ^{n } v_{h_n}(j,k)\\
 &=(n-1)^{-1}   \sum _ { j =1 } ^{n }  v_{h_n}(j,j)
 + (n-1)^{-2}   (n-2) \sum _ { j =   1} ^{n } \sum _ {  k\neq  j } ^{n } v_{h_n}(j,k).
\end{align*}
The term associated with double summation over $j$ and $k$ is a degenerate $U$-statistic, as $\mathbb E[ v_{h_n}(j,k)\mid X_k ] =\mathbb E[ v_{h_n}(j,k)\mid X_j ]=0$. Consequently, following the same lines as in the treatment of the first term of $A_n $, and using (\ref{ineq:order2moments_bis_2}) in Lemma \ref{lemma:auxiliary_bound_kernel}, we get
\begin{align*}
\mathbb E\left[\left(  (n-1)^{-2}   (n-2) \sum _ { j =   1} ^{n } \sum _ {  k\neq  j } ^{n } v_{h_n}(j,k)  \right)^2  \right] = O(1) \mathbb E [v_{h_n}(1,2)^2 ] =   O( h_n^{-d})
\end{align*}
which goes to $0$, when multiplied by $h_n^{2d}$. The remaining term is a sum of independent and identically distributed random variables. We have, by computing the variance of the centred average,
 \begin{align*}
(n-1)^{-1}    \sum _ { j =1 } ^{n }  (v_{h_n}(j,j)  =  O_{\mathbb P} (n^{-1/2}) \sqrt { \mathbb E v_{h_n}(1,1)^2 }  + n (n-1)^{-1}  \mathbb E [ v_{h_n} (1,1)],
 \end{align*}
where the first term, using (\ref{ineq:order2moments_bis_0}) in Lemma \ref{lemma:auxiliary_bound_kernel}, is $O(n^{-1/2}h_n^{-d})$ which goes to $0$ when multiplied by $h_n$. The dominating term is in fact the last one, as by Lemma \ref{lemma:auxiliary_limit_kernel}, it holds that $ h_n ^d \mathbb E [ v_{h_n}(1,1)]\to v_K\lambda(S)$.

\paragraph{Proof that $ h_n^d C_n \to 0$ in probability.}
We are based on similar decompositions as for $A_n$ involving $U$-statistics. Let $\ell_{h_n}(x) = (f_{h_n}(x)-f_{0}(x)+\Delta_nq(x))/f_{0}(x)$ and note that in virtue of Lemma \ref{lemma:auxiiliary_f_conv}, it holds
\begin{align*}
\|\ell_{h_n}\|_S \leq b^{-1}  \left( h_n^2 \| g\|_{\mathbb R^d} \int \|u\|_2^2 K(u) \,\diff u +\Delta_n \|q\|_{\mathbb R^d} \right).
\end{align*}
Then
\begin{align*}
C_n =   (n-1)^{-1}\sum _ { i = 1 } ^{n } \sum _ { j \neq  i  } ^{n } (a_{h_n}(i,j)\ell_{h_n,i} \mathrm 1_{\{ X_i\in S\}} - b_{h_n}(j )  )+  \sum _ { i = 1 } ^{n } b_{h_n}(i )  ,
\end{align*}
with $b_{h_n}(j ) =  \mathbb E[ a_{h_n}(i,j)\ell_{{h_n},i}  \mathrm 1_{\{ X_i\in S\}}  \mid X_j] $.
The term on the left is a degenerate $U$-statistic for which it holds
\begin{align*}
\mathbb E\left[\left(  (n-1)^{-1}\sum _ { i = 1 } ^{n } \sum _ { j \neq  i  } ^{n } (a_{h_n}(i,j)\ell_{h_n,i} \mathrm 1_{\{ X_i\in S\}} -b_{h_n}(j ) )\right)^2  \right] = O( 1)\mathbb E[a_{h_n}(1,2)^2 \ell_{{h_n},1}^2 \mathrm 1_{\{X_1 \in  S \}} ].
\end{align*}
Using (\ref{ineq:order2moments00}) in Lemma \ref{lemma:auxiliary_bound_kernel} and the previous bound for $\|\ell_{h_n}\|_S$, we find
\begin{align*}
\mathbb E[a_{h_n}(1,2)^2 \ell_{{h_n},1}^2 \mathrm 1_{\{X_1 \in  S \}} ] & =  \mathbb  E\left [\left(\frac{ V_{h_n}(X_1) \ell_{{h_n},1}^2}{f_{0,1} ^2  } \right)\mathrm 1_{\{X_1 \in  S \}} \right] = O(h_n^{-d}(h_n^2 + \Delta_n)) ,
\end{align*}
where $V_h$ is defined in (\ref{def:V_h}). The previous bound multiplied by $h_n^{2d}$ goes to $0$. Using that $\mathbb E[a_{h_n}(1,2)\ell_{h_n,1} \mathrm 1_{\{ X_1\in S\}}]  = 0$ and (\ref{ineq:order2moments2}), the variance of the term on the right in $C_n$ is smaller than 
\begin{align*}
n\mathbb E[\mathbb E[ a_{h_n}(1,2)\ell_{h_n,1} \mathrm 1_{\{ X_1\in S\}} \mid X_2] ^2 ] &\leq n\mathbb E[\mathbb E[ |a_{h_n}(1,2)| \mid X_2 ] ^2 ]\|\ell_{h_n}\|_S ^2   
 = O(n (h_n^4 + \Delta_n^2))
\end{align*}
which, multiplied by $h_n^{2d}$, goes to $0$ by hypothesis. Hence $h_n^{d}C_n\to 0$, in probability and the proof is complete.
\end{proof}

\begin{app_lemma}\label{lemma:cvlin}
Under assumptions (H\ref{ass:base_density}), (H\ref{ass:kernel}) and (H\ref{ass:Hcvlin_f_0_St}), 
if $  n h_n^d  \Delta_n \to 0$, we have
\begin{align*}
&  {h_n^{d}} \, \sum _ { i = 1 } ^{n }  \left(\frac{ \hat f _{n,i} ^{\text{LR}} - f_{0,i}}{f_{0,i}} \right) = o_{\mathbb P} (1). 
\end{align*} 
\end{app_lemma}

\begin{proof}
The decomposition is as follows
\begin{align}
\nonumber & h_n^{d}\sum _ { i = 1 } ^{n }  \left(\frac{ \hat f _{n,i} ^{\text{LR}} - f_{0,i}}{f_{0,i}} \right)  \\
 \label{eq:decomp_lemma_lin_np} &=  h_n^{d} (n-1) ^{-1} \sum _ { i = 1 } ^{n } \sum _ { j\neq  i } ^{n } \left(\frac{ K_{h_n} (i,j) - f_{h_n,i } }{f_{0,i}} \right)  +  h_n^{d}\sum _ { i = 1 } ^{n }\left(\frac{ f_{h_n,i }-f_{0,i} }{f_{0,i}} \right) +  h_n^{d}\Delta_n \sum _ { i = 1 } ^{n }\frac{  q_i}{f_{0,i}} .
\end{align}
The expectation of the last term is $n h_n^{d}\Delta_n \int  q(x) \,\diff x$ which goes to $0$ by assumption. We can now focus on the first and second term of the decomposition.

\noindent\textit{Treatment of the second term in (\ref{eq:decomp_lemma_lin_np}).}
Using that $\int K(u) \,\diff u = 1$, the considered term is a centred empirical sum. Using Lemma \ref{lemma:auxiiliary_f_conv}, its variance is then bounded by
\begin{align*}
\mathbb E \left[\left( h_n^{d}\sum _ { i = 1 } ^{n }\frac{ f_{h_n,i }-f_{0,i} }{f_{0,i}}  \right) ^2 \right] 
&\leq  nh_n^{2d} \int  \frac{(f_{h_n }(x)-f_{0}(x))^2}{f_0(x)}\,\diff x  \\
&\leq nh_n^{2d+4} \int  \frac{g(x)^2 }{f_0(x)} \,\diff x  \left( \int u^2    K(u)  \,\diff u\right)^2,
\end{align*}
which goes to $0$.

\noindent\textit{Treatment of the first term in (\ref{eq:decomp_lemma_lin_np}).}
Using that $\int K(u) \,\diff u = 1$, one can verify that it is a degenerate $U$-statistic.
Here the variance can not be computed directly because the leading term $\mathbb E \left[ \frac{ (K_{h_n} (1,2) - f_{h_n,1 })^2 }{f_{0,1}^2}   \right] $ is not necessarily finite. Hence we decompose according to the $X_i$ in $S_{b_n} $ and the others, with $b_n = (\epsilon /  nh_n^d)^{1/\beta}$ where $\beta $ is given in (H\ref{ass:Hcvlin_f_0_St}) and $\epsilon > 0$. We introduce  
\begin{align*}
k(x,y) = \frac{ K_{h_n} (x,y)  }{f_{0}(x)},
\end{align*}
and define the linear operator $Q_{\mathbb P}:L_2(\mathbb P) \to L_2(\mathbb P)$ as 
\begin{align*}
Q_{\mathbb P}[w] (x,y) = w(x,y) - \mathbb E[w(x,X_1)] - \mathbb E[ w(X_1,y)] + \mathbb E[w(X_1,X_2)].
\end{align*}
Because $\mathbb E[ k(X_1,y)] = \mathbb E[k(X_1,X_2)]=1$ for all $y\in\mathbb R^d$, one sees that
\begin{align*}
  \sum _ { i = 1 } ^{n } \sum _ { j\neq  i } ^{n } \left(\frac{ K_{h_n} (i,j) - f_{h_n,i } }{f_{0,i}} \right) 
&= \sum _ { i = 1 } ^{n } \sum _ { j\neq  i } ^{n }  Q_{\mathbb P} (k)_{i,j}  \\
& =   \sum _ { i = 1 } ^{n } \sum _ { j\neq  i } ^{n } \left(  Q_{\mathbb P}  (k \mathrm 1_{S_{b_n}})_{i,j}  +   Q_{\mathbb P} (k \mathrm 1_{S_{b_n}^c})_{i,j}  \right) .
\end{align*}
Because the summation over $Q_{\mathbb P}  (k \mathrm 1_{S_{b_n}})$ is a degenerate $U$-statistics, we get that
\begin{align*}
\mathbb E \left[ \left( h_n^{d} (n-1)^{-1}  \sum _ { i = 1 } ^{n } \sum _ { j\neq  i } ^{n }   Q_{\mathbb P}  (k \mathrm 1_{S_{b_n}})_{i,j} \right)^2\right] 
&= O(h_n^{2d})\mathbb E \left[ \frac{ (K_{h_n} (1,2) - f_{h_n,1 })^2 }{f_{0,1}^2}1_{S_{b_n}}(X_1)   \right].
\end{align*}
Defining the kernel $\tilde K = K^2/v_K$ and $\tilde f_h =  f_0\star \tilde K_h$, we obtain
\begin{align*}
\mathbb E \left[ \frac{ (K_{h_n} (1,2) - f_{h_n,1 })^2 }{f_{0,1}^2}1_{S_{b_n}}(X_1)   \right]  &\leq \mathbb E \left[ \frac{\mathbb E[ K_{h_n} (1,2)^2\mid X_1] }{f_{0,1}^2}1_{S_{b_n}}(X_1)   \right] \\
&= v_Kh_n^{-d} \int_{S_{b_n}} \frac{\tilde f_{h_n}(x) }{f_{0}(x) } \,\diff x \\
& = v_Kh_n^{-d} \left( \int_{S_{b_n}} \int \frac{f(x-h_nu)}{f_{0}(x) } \tilde K (u) \diff u \,\diff x\right)\\
& \leq   v_Kh_n^{-d} \lambda(S_{b_n}) \left(  \sup_{x\in S_{b_n}} \sup_{u\in[-1,1]^d } \frac{f(x+h_nu)}{f_{0}(x) } \right).
\end{align*}
For the term with $Q_{\mathbb P} (k \mathrm 1_{S_{b_n}^c}) $, we obtain that
\begin{align*}
\mathbb E\left[\left| \sum _ { i = 1 } ^{n } \sum _ { j\neq  i } ^{n }   Q_{\mathbb P} (k \mathrm 1_{S_{b_n}^c})_{i,j} \right|  \right]
&\leq n(n-1)\mathbb E [ |Q_{\mathbb P} (k \mathrm 1_{S_{b_n}^c})_{1,2}| ]\\
&\leq 4 n(n-1)\mathbb E [ | k_{1,2}|1_{S_{b_n}^c}(X_1) ]\\
&  =  4 n(n-1)\int_{S_{b_n}^c}  f_h(x) \,\diff x.
\end{align*}
From Lemma \ref{lemma:auxiliary_int_f_h}, we deduce that $\int_{S_{b_n}^c}  f_h(x) \,\diff x  \leq c_2 b_n ^\beta  = c_2 \epsilon / nh_n^d$  . To conclude, we have shown that there exists a constant $\tilde C>0$ such that
\begin{align*}
&\mathbb E\left[\left| h_n^d (n-1)^{-1} \sum _ { i = 1 } ^{n } \sum _ { j\neq  i } ^{n }   Q_{\mathbb P}  (k )\right| \right] \leq \tilde C\left( \sqrt{ h_n^{d } \lambda(S_{b_n}) } +n h_n^d b_n^\beta    \right)  =  \tilde C\left( \sqrt{ h_n^{d } \lambda(S_{b_n}) } +\epsilon  \right) .
\end{align*} 
Invoking (H\ref{ass:Hcvlin_f_0_St}) and because $\epsilon$ is arbitrarily small, the limit as $n\to\infty$ is $0$.
\end{proof}

\begin{app_lemma}\label{lemma:cvquad_parametric}
Under (A\ref{ass:parametric_model_consistency}) and (A\ref{ass:parametric_model_assymptotic_normality}), we have
\begin{align*}
&\sum _ { i = 1 } ^{n }  \left(\frac{ f_{\hat \theta_n,i} - f_{0,i}}{f_{0,i}} \right)^2\mathrm 1_{\{X_i\in S\}}  = O_{\mathbb P} (1) .
\end{align*} 
\end{app_lemma}
\begin{proof}
Using (\ref{eq:lip_theta_f}), we have that, with probability going to $1$,
\begin{align*}
 \sum _ { i = 1 } ^{n }  \left(\frac{ f_{\hat \theta_n,i} - f_{0,i} }{f_{0,i}} \right)^2 \mathrm 1_{\{X_i\in S\}} &\leq   \left(n^{-1}\sum _ { i = 1 } ^{n }\frac{\dot \ell (X_i)^2  \sup_{\theta\in B(\theta_0,\delta)}  f_\theta (X_i)^2}{f_{0,i}^2}\mathrm 1_{\{X_i\in S\}}\right)  n\|\hat \theta_n- \theta_0\|_2^2\\
 &\leq \|\dot \ell   \sup_{\theta\in B(\theta_0,\delta)}  f_\theta  \|^2_{\mathbb R^d}b^{-2}n\|\hat \theta_n- \theta_0\|_2^2
\end{align*}
which is a tight sequence in light of (\ref{eq:rep_theta}).
\end{proof}

\begin{app_lemma}\label{lemma:cvlin_param}
Under (A\ref{ass:parametric_model_consistency}) and (A\ref{ass:parametric_model_assymptotic_normality}), we have
\begin{align*}
& \sum _ { i = 1 } ^{n }  \left(\frac{ f_{\hat \theta_n,i} - f_{0,i} }{f_{0,i}} \right) = O_{\mathbb P}(1).
\end{align*}
\end{app_lemma}

\begin{proof}

In virtue of (A\ref{ass:parametric_model_assymptotic_normality}), the map $\theta\mapsto \log f_\theta(x)$ is differentiable at $\theta_0$, for $P$-almost every $x\in\mathbb R^d $ with derivative $\dot \ell_{\theta_0}(x)$ (this is obtained in \cite{vandervaart:1998}, in the proof of Theorem 5.39). Using stability properties for the composition, the map $\theta\mapsto  f_\theta(x) = \exp(\log(f_\theta(x) ) )$ is differentiable at $\theta_0$, for $P$-almost every $x\in\mathbb R^d$ with derivative $\dot \ell_{\theta_0}(x) f_0(x)$. We are in position to apply Lemma 19.31 in \cite{vandervaart:1998}, with $r_n = \sqrt n $ and $m_\theta = f_\theta/f_0$. From the mentioned lemma, as $\sqrt n (\hat \theta_n - \theta_0)$ is tight, defining
\begin{align*}
&T_i(\theta) = \left[ \left(\frac{ f_{\theta,i} - f_{0,i} }{f_{0,i}} \right)  - (\theta-\theta_0) ^T  {\dot \ell_{\theta_0,i}}\right],\\
&t(\theta) = \int \left[{ f_{\theta}(x) - f_{0}(x) }   - (\theta -\theta_0) ^T {\dot \ell_{\theta_0}(x) f_0(x) }\right] \,\diff x,
\end{align*}
we obtain
\begin{align*}
\left|  \sum _ { i = 1 } ^{n } \{ T_i(\hat \theta_n)   - t(\hat \theta_n) \}  \right| =o_{\mathbb P}(1). 
\end{align*}
Actually, recalling that $\int \dot \ell _{ \theta_0} (x) f_0(x)\,\diff x = 0 $, we find that, for all $\theta\in \Theta $, $t(\theta) = 0$. Hence, we obtain
\begin{align*}
 \sum _ { i = 1 } ^{n }  \left(\frac{ f_{\hat \theta_n,i} - f_{0,i} }{f_{0,i}} \right) = n^{1/2} (\hat \theta_n-\theta_0) ^T\left[ n^{-1/2} \sum _ { i = 1 } ^{n }  { \dot \ell _{ \theta_0, i} }  \right]+o_{\mathbb P}(1).
\end{align*}
where the first term is a $O_{\mathbb P} (1)$. 
\end{proof}

\subsection{Auxiliary results}\label{app:analytical_results}

Recall some definitions, for any $h>0$,
\begin{align}
\label{def:V_h}&V_h(X_1) =  \mathbb E [ (K_{h}(1,2) - f_{h,1})^2 \mid X_1], \\
\nonumber&a_h(x,y) = \frac{ K_h(x-y)  - f_{h}(x)}{f_{0}(x)  } ,\\
\nonumber &u_{h}(x,y,z) = a_h(x,y) a_h(x,z) \mathrm 1_{\{x \in  S \}},
\end{align}
as well as the short-cut $g(i,j,k) $ for $g(X_i,X_j,X_k) $.

\begin{app_lemma}\label{lemma:auxiliary_bound_kernel}
Under (H\ref{ass:base_density}) and (H\ref{ass:kernel}), if $S\subset \mathbb R^d $ is such that for all $x\in S$, $f_0(x)>b>0$, we have, for any $h>0$,
 \begin{align}
 \label{ineq:order2moments00}
 &V_h(X_1) \leq h^d C_0 \\
 \label{ineq:order2moments0} 
& \mathbb E \left[ \mathbb E [ u _h(1,2,2) \mid X_1]^2 \right]  \leq h^{-2d} C_1 \\
\label{ineq:order2moments1} 
&\mathbb E[{u}_{h}(1,2,3)^2 ]  \leq h^{-2d} C_1 \\
\label{ineq:order2moments2} 
& \mathbb E[ |a_h(1,2)| \mid X_2]  \leq  2  \\
 \label{ineq:order2moments_bis_0}
&\mathbb E [ \, \mathbb E[u_{h}(1,2,2) \mid X_2 ] ^2\, ] \leq h^{-2d} C_2 + C_3\\
\label{ineq:order2moments_bis_1} 
&\mathbb E [\, u_{h}(1,2,2)^2 \, ]   \leq  h^{-3d} C_4 +C_5\\
\label{ineq:order2moments_bis_2} 
& \mathbb E \left[  \, \mathbb E[u_{h}(1,2,3) \mid X_2,X_3 ] ^2\,   \right] \leq  h^{-d} C_6 + C_6
\end{align}
where the constants $C_k$, $k=0,\ldots, 7$ depends on $K$ and $f_0$ only.
\end{app_lemma}

\begin{proof}
Remark that because $K$ is bounded and $\int |K(u)|\,\diff u<\infty$, we have $\int |K(u)|^k\,\diff u<\infty$, for any $k \geq 1$. Note that, for every $h>0$,
\begin{align*}
V_h(X_1) \leq\mathbb E [ K_{h}(2,1) ^2\mid X_1] \leq   h^{-d}v_K  \|f_0\|_{\mathbb R^d}.
\end{align*}
We obtain (\ref{ineq:order2moments0}) by writing
\begin{align*}
\mathbb E \left[ \mathbb E [ u _h(1,2,2)\mid  X_1]^2 \right]  =   \mathbb E\left[ \frac{V_h(X_1)^2}{f_{0,1}^4} \mathrm 1_{\{X_1 \in  S \}} \right] \leq  h^{-2d}v_K^2  \|f_0\|_{\mathbb R^d}^2 b^{-4}.
\end{align*}
To establish (\ref{ineq:order2moments1}), note that
\begin{align*}
&\mathbb E[ u_h(1,2,3) ^2] = \mathbb E\left[ \frac{V_h(X_1)^2}{f_{0,1}^4} \mathrm 1_{\{X_1 \in  S \}} \right] = \mathbb E \left[ \mathbb E [ u _h(1,2,2)| X_1 ]^2 \right]  .
\end{align*}
For (\ref{ineq:order2moments2}), write
 \begin{align*}
 \mathbb E[ |a_h(1,2)| \mid X_2] &=    \int { |  K_{h} (x-X_2 ) - f_{h}(x)| } \,\diff x   \leq   \int     K_{h} (x-X_2 ) \,\diff x     +   \int f_{h}(x) \,\diff x  = 2.
\end{align*}
Inequality (\ref{ineq:order2moments_bis_0}) follows from the lines
\begin{align*}
\mathbb E [ \, \mathbb E[u_{h}(1,2,2) \mid X_2 ] ^2\, ]&= \int  \left( \int   \frac{ (K_{h}(x-y) - f_{h}(x) )^2}{f_{0}(x)} \mathrm 1_{\{x \in  S \}}  \,\diff x  \right)^2 f_0(y) \,\diff y\\
&\leq 2 \int  \left( \int   \frac{ K_{h}(x-y)^2 + f_{h}(x) ^2}{f_{0}(x)} \mathrm 1_{\{x \in  S \}}  \,\diff x  \right)^2 f_0(y) \,\diff y\\
&\leq 2b^{-2}  \int  \left( \int   { K_{h}(x-y)^2 + f_{h}(x) ^2}   \,\diff x  \right)^2 f_0(y) \,\diff y\\
&\leq 2b^{-2}   \int  \left( h^{-d} v_K  + \|f_h\|_{\mathbb R^d} \right)^2    f_0(y) \,\diff y\\
&  \leq 4b^{-2}    \left( h^{-2d} v_K^2   + \|f_0\|_{\mathbb R^d}^2  \right) .
\end{align*}
To show (\ref{ineq:order2moments_bis_1}), write
\begin{align*}
&\mathbb E[ u_h(1,2,2) ^2] = \mathbb E\left[ \frac{(K_{h}(1,2) - f_{h,1})^4}{f_{0,1}^4} \mathrm 1_{\{X_1 \in  S \}} \right].
\end{align*}
Using that $(a+b)^4\leq 8(a^4+b^4)$, we obtain
\begin{align*}
\nonumber \mathbb E [u_{h}(1,2,2)^2 ]   &\leq 8\mathbb E \left[ \left(\frac{ K_{h}(1,2) ^4 + f_{h,1} ^4 }{f_{0,1}^4} \right) \mathrm 1_{\{X_1\in  S \}}\right] \\
&\leq 8 b^{-4} \left( h^{-3d} \|f_0\|_{\mathbb R^d} \int K(u)^4\,\diff u +\|f_0\|_{\mathbb R^d}^4 \right).
\end{align*}
For (\ref{ineq:order2moments_bis_2}), we have
\begin{align*}
&\mathbb E \left[  \, \mathbb E[u_{h}(1,2,3) \mid X_2,X_3 ] ^2\,   \right] \\
&=  \mathbb E \left[   \mathbb E\left[  \left(\frac{ K_{h}(3,1) - f_{h,3}}{f_{0,3}} \right) \left(\frac{ K_{h}(3,2) - f_{h,3}}{f_{0,3}} \right)\mathrm 1_{\{X_3 \in  S \}} \mid X_1,X_2  \right] ^2\right].
\end{align*}
We develop and compute bounds for each term. The larger term will be the one associated with the product of the kernels. We have, by Jensen's inequality, for any $(y,z)\in \mathbb R^d\times \mathbb R^d$, 
\begin{align*}
\psi_h(y,z) &:= \left(\int  \left(\frac{ K_{h}(x-y)  K_{h}(x-z) }{f_{0}(x)} \right)\mathrm 1_{\{x \in  S \}}\,\diff x\right)^2\\ 
  &\leq b^{-2} \left(\int  { K_{h}(x-y)  K_{h}(x-z) } \,\diff x \right)^2\\ 
 &= b^{-2} h^{-2d}  \left( \int   K(u)  K((y-z)/h+u) \,\diff u \right)^2\\
 & \leq b^{-2} h^{-2d}   \int   K(u)  K((y-z)/h+u)^2 \,\diff u .
 \end{align*}
Then we obtain
\begin{align*}
&\mathbb E \left[  \mathbb E\left[  \left(\frac{ K_{h}(3,1) }{f_{0,3}} \right) \left(\frac{ K_{h}(3,2) }{f_{0,3}} \right)\mathrm 1_{\{X_3 \in  S \}} \mid X_1,X_2  \right] ^2 \right ] \\
&=\int \int \psi_h(y,z) f_0(y)f_0(z) \,\diff y\diff z\\
&\leq   b^{-2} h^{-2d}  \int \int \int  { K(u)   K((y-z)/h+u) ^2} f_0(y)f_0(z) \,\diff y\diff z\diff u\\
&=  h^{-d} b^{-2}  \int \int \int  { K(u)  K(v+u)^2} f_0(z+hv)f_0(z) \,\diff v\diff z\diff u\\
&\leq  h^{-d} b^{-2} \|f_0\|_{\mathbb R^d}  \int \int   { K(u)  K(v+u) ^2}  \,\diff v \diff u\\
& =  h^{-d} b^{-2} \|f_0\|_{\mathbb R^d}   v_K .
\end{align*}
Moreover, as
\begin{align*}
 \mathbb E\left[  \left(\frac{ K_{h}(3,1) }{f_{0,3}} \right) \left(\frac{ f_{h,3} }{f_{0,3}} \right)\mathrm 1_{\{X_3 \in  S \}} \mid X_1,X_2  \right] & =\int  \left(\frac{ K_{h}(x-X_1)  f_{h} (x) }{f_{0}(x)} \right)\mathrm 1_{\{x \in  S \}}\,\diff x
  \leq  \|f_0\|_{\mathbb R^d} b^{-1},
 \end{align*}
 and 
\begin{align*}
 \mathbb E\left[  \left(\frac{ f_{h,3}}{f_{0,3}} \right)^2 \mathrm 1_{\{X_3 \in  S \}} \mid X_1,X_2  \right] &= \int  \frac{ f_{h} (x)^2 }{f_{0}(x)} \mathrm 1_{\{x \in  S \}}\,\diff x 
 \leq \|f_0\|_{\mathbb R^d}   b^{-1}  .
 \end{align*}
 we finally obtain the result. 
\end{proof}

\begin{app_lemma}\label{lemma:auxiliary_limit_kernel}
Under (H\ref{ass:base_density}) and (H\ref{ass:kernel}), if $S\subset \mathbb R^d$ is such that for all $x\in S$, $f_0(x)>b>0$, we have that
\begin{align*}
\lim _{ h\rightarrow 0} h^d \mathbb E \left[ \left( \frac{ K_{h}(1,2) - f_{h,1}}{f_{0,1}} \right)^2\mathrm 1_{\{X_1 \in  S \}}\right]    = v_K \lambda(S).
\end{align*}

\end{app_lemma}

\begin{proof}
Write
\begin{align*}
\mathbb E \left[ \left( \frac{ K_{h}(1,2) - f_{h,1}}{f_{0,1}} \right)^2\mathrm 1_{\{X_1 \in  S \}}\right]
& = \mathbb E \left[ \frac{ V_h(X_1)}{f_{0,1}^2} \mathrm 1_{\{X_1 \in  S \}}\right]\\
& = \mathbb E \left[  \frac{ \mathbb E [ K_{h}(1,2)^2 \mid X_1] }{f_{0,1}^2} \mathrm 1_{\{X_1 \in  S \}} \right] - \mathbb E \left[  \frac{f_{h,1}^2}{f_{0,1}^2} \mathrm 1_{\{X_1 \in  S \}} \right].
\end{align*}
The right-hand side is bounded by $b^{-2}\|f_0 \|_{\mathbb R^d} $, hence its participation in the stated limit is $0$. For the left-hand side term, use $\tilde K = K^2/v_K$ and write
\begin{align*}
h^d \mathbb E \left[  \frac{\mathbb E [ K_{h}(1,2)^2 \mid X_1] }{f_{0,1}^2} \mathrm 1_{\{X_1 \in  S \}}  \right]  &= h^d  \int\int  \frac{ f_0(y) }{f_{0}(x)}K_{h}(x-y)^2 \mathrm 1_{\{x \in  S \}} 
\,\diff y \diff x\\
&=  \int\int  \frac{f_0(x-hu)}{f_{0}(x)}  K(u)^2 \mathrm 1_{\{x \in  S \}} 
\,\diff u  \diff x\\
&=  v_K \int\int   \frac{ f_0(x-hu) }{f_{0}(x)}\tilde K(u) \mathrm 1_{\{x \in  S \}} 
\,\diff u \diff x\\
& =  v_K \lambda(S) +  v_K \int \int  \frac{  (f_0(x-hu)-f_0(x)) }{f_{0}(x)} \tilde K(u) \mathrm 1_{\{x \in  S \}} 
\,\diff u \diff x.
\end{align*}
It remains to note that the term in the right goes to $0$, as $h\to 0$, in virtue of the Lebesgue dominated convergence theorem. 
\end{proof}

\begin{app_lemma}\label{lemma:auxiiliary_f_conv}
Under (H\ref{ass:base_density}) and (H\ref{ass:kernel}), we have, for every $x\in \mathbb R^d$ and $h>0$,
\begin{align*}
\nonumber |f_0 \star K_h (x) - f_0(x) | &\leq  g(x) h^2 \int \|u\|_2^2   K(u)  \,\diff u .
\end{align*}
\end{app_lemma}

\begin{proof}
Note that $\int K(u)\, \diff u =1$ and, by symmetry,  $\int u K(u)\,\diff u = 0$. Write
\begin{align*}
 |f_0 \star K_h (x) - f_0(x) | &=\left|  \int (f_0(x-hu) - f_0(x)) K(u) \,\diff u\right|\\
&= \left|  \int (f_0(x-hu) - f_0(x) -(hu)^T \nabla f_0(x) ) K(u) \,\diff u\right|\\
&\leq   \int | f_0(x-hu) - f_0(x) -(hu)^T \nabla f_0(x)| \,   K(u)  \,\diff u ,
\end{align*}
and use (H\ref{ass:base_density}) to conclude.
\end{proof}

\begin{app_lemma}\label{lemma:auxiliary_int_f_h}
Under (H\ref{ass:base_density}), (H\ref{ass:kernel}) and (H\ref{ass:Hcvlin_f_0_St}), there exists $c_2>0$ such that $\int_{S_{b_n}^c}  f_{h_n}(x) \,\diff x  \leq c_2 b_n^\beta $.
\end{app_lemma}

\begin{proof}
Note that
\begin{align*}
\int_{S_{b_n}^c}  f_h(x) \,\diff x &= \mathbb P (S_{b_n}^c) +  \int_{S_{b_n}^c} ( K_{h_n} \star f_0 -f_0)  \,\diff x\\
& =   \mathbb P (S_{b_n}^c) +  \int_{} \, ( K_{h_n} \star \mathrm 1_{S_{b_n}^c}(x) -\mathrm 1_{S_{b_n}^c}(x) ) f_0(x) \,\diff x.
\end{align*}
The term in the left is bounded by $cb_n^{\beta}$ as supposed in (H\ref{ass:Hcvlin_f_0_St}).
For the term in the right, define $S_{b_n,h_n} = \{y+h_nu\,:\, u\in [-1,1]^d,\, y\in S_{b_n}\}$. Note that, by (H\ref{ass:kernel}), as soon as $x\notin S_{b_n,h_n}$, $ K_{h_n} \star \mathrm 1_{S_{b_n}^c}  (x)  = 1$, hence
\begin{align*}
|  K_{h_n} \star \mathrm 1_{S_{b_n}^c}  (x) -\mathrm 1_{S_{b_n}^c}(x)  |\leq  \mathrm 1_{S_{b_n,h_n}}   . 
\end{align*}
Moreover, for any $x  \in S_{b_n,h_n} $, we have, by (H\ref{ass:Hcvlin_f_0_St}), that
\begin{align*}
f_0(x) &\leq   \sup_{y\in  S_{b_n}} \sup_{u\in[-1,1]^d}  f_0(y+h_nu) \leq  b_n \sup_{y\in  S_{b_n}} \sup_{u\in[-1,1]^d} \frac{ f_0(y+h_nu) }{f_0(y)}  = b_n C,
\end{align*}
hence, $\mathrm 1_{S_{b_n,h_n}}  \leq \mathrm 1_{f_0(x)\leq C b_n } $, leading to
\begin{align*}
\left| \int_{}\, ( K_{h_n} \star \mathrm 1_{S_{b_n}^c}(x) -\mathrm 1_{S_{b_n}^c}(x) ) f_0(x) \,\diff x \right| \leq  \int  \mathrm 1_{f_0(x)\leq C b_n }   f_0(x)\,\diff x \leq (C b_n) ^\beta . 
\end{align*}

\end{proof}

\section{Parametric maximum likelihood estimator}\label{app:parametric}

In this section are reported some classical results on the maximum likelihood estimator of the density. When the model is well-specified, we need the consistency and the asymptotic normality of the estimated parameter $\theta_0$.

\begin{enumerate}[(\text{A}1)]
\item \label{ass:parametric_model_consistency} The set $\Theta\subset \mathbb R^q$ is compact. The model $\mathcal P= \{ f_\theta\ : \ \theta \in \Theta \}$, a collection of densities on $\mathbb R^d$, is identifiable, i.e., for every $\theta_1\neq \theta_2 $ in $ \Theta$, $f_{\theta_1} \neq  f_{\theta_2}$ and the envelop $F_\Theta (x)= \sup_{\theta\in \Theta} f_\theta(x)$ is such that $\mathbb E [\log(F_{\Theta,1})]<+\infty$. There exists an $\mathbb R^+$-valued measurable function $\dot{\ell}$ with $E \dot{\ell} (X_1)^2 <\infty$ for every $x\in \mathbb R^d$, for every $\theta_1$ and $\theta_2$ in $\Theta$,
\begin{align*}
|\log(f_{\theta_1} (x))-\log(f_{\theta_2} (x))| \leq \dot \ell (x) \| \theta_1-\theta_2\|_2.
\end{align*}
There exists $\delta>0$ such that the function $\dot \ell \times  \sup_{\theta\in B(\theta_0,\delta)}  f_\theta $ is bounded.
\setcounter{ass:A}{\value{enumi}}
\end{enumerate}

It follows from (A\ref{ass:parametric_model_consistency}) that the class of functions $\mathcal P$ is Glivenko-Cantell \citep[Theorem 2.7.11]{wellner1996}, i.e.,
\begin{align}\label{eq:GC}
\sup_{\theta \in \Theta} \left| n^{-1}\sum_{i=1} ^n (\log(f_\theta(X_i) - E[\log(f_\theta(X_1)  ])\right| \rightarrow 0.
\end{align}
Whenever $ f_0\in \mathcal P$, it holds that $\hat \theta_n\rightarrow \theta_0$, in probability \citep[Theorem 2.1]{newey+m:1994} or \citep[Lemma 5.35]{vandervaart:1998}. For now, asking the above Lipschitz condition to guarantee the Glivenko-Cantelli might seem a bit restrictive \citep[Lemma 2.4]{newey+m:1994}, but this condition will also be required to derive asymptotic normality of $\hat \theta_n$ as well as to obtain uniform convergence (over $x\in \mathbb R^d$) of $f_{\hat \theta_n}(x)$ to $f_{\theta_0}(x)$. Indeed, we have that for any $\delta>0$, with probability going to $1$, $\hat \theta_n \in B(\theta_0,\delta)$. Hence, using the mean-value theorem, we find 
\begin{align}\label{eq:lip_theta_f}
|f_{\hat \theta_n}(x) - f_{\theta_0}(x)| \leq   \|\hat \theta_n -\theta_0\|_2 \dot \ell (x) \sup_{\theta\in B(\theta_0,\delta)}  f_\theta (x)
\end{align}
for every $x\in \mathbb R^d$. Conclude using that $\dot \ell \times  \sup_{\theta\in B(\theta_0,\delta)}  f_\theta $ is bounded and the convergence in probability of $\hat \theta_n$ to $\theta_0$.

\begin{enumerate}[(\text{A}1)]
\setcounter{enumi}{\value{ass:A}}
\item \label{ass:parametric_model_assymptotic_normality}
The true parameter $\theta_0$ an interior point of $\Theta\subset \mathbb R^q$. The model $\mathcal P$ is differentiable in quadratic mean at $\theta_0$, i.e., there exists a measurable vector-valued function $\dot \ell_{\theta_0}$, with $ E[ \| \dot \ell_{\theta_0}(X_1)\|_2^2]$, such that
\begin{align*}
\int  \left[ \sqrt { f_{\theta}}- \sqrt { f_{\theta_0}}-\frac 1 2 (\theta-\theta_0) ^T \dot \ell_{\theta_0} \sqrt{f_{\theta_0}}  \right]^2 \diff \lambda   = o(\|\theta-\theta_0\|_2^2).
\end{align*}
The matrix $  \mathcal I = E[  \dot \ell_{\theta_0}(X_1) \dot \ell_{\theta_0}(X_1)^T]$ is invertible. 
\end{enumerate}

As a consequence of the previous set of conditions \cite[Lemma 5.39]{vandervaart:1998}, we have
\begin{align}\label{eq:rep_theta}
n^{1/2}(\hat \theta_n- \theta_0 )= \mathcal I^{-1}  n^{-1/2} \sum _ { i = 1 } ^{n } \dot \ell_{\theta_0}(X_i) + o_{\mathbb P} (1) .
\end{align} 
where $E[\dot \ell_{\theta_0}(X_1) ] = 0$. In particular, it holds that $\sqrt n (\hat \theta_n- \theta_0 ) = O_{\mathbb P}(1)$.

\end{appendices}

\paragraph{Acknowledgments.} 
The authors are grateful to Ingrid Van Keilegom, Anouar El Ghouch and Sylvie Huet for useful comments and references. The authors are also grateful to Christian Robert for some motivating remarks at the very beginning of this work.

\bibliographystyle{abbrv} 
\bibliography{bib_mixture}

\end{document}